\documentclass[preprint]{elsarticle}
\usepackage{cmap} 
           \usepackage[T2A]{fontenc}   
           \usepackage[utf8]{inputenc}
           \usepackage[english]{babel}
\usepackage{amsmath}
\usepackage{amsfonts}
\usepackage{amssymb}
\usepackage{amsthm}
\usepackage[pdftex,unicode]{hyperref}
\usepackage{amscd}
\usepackage[all,cmtip]{xy}
\usepackage{mathtools}
\usepackage{comment}
\usepackage{tikz-cd}
\usetikzlibrary{babel}

  \def\R{\mathbb R}
\def\T{\mathbb T}  
\def\C{\mathbb C}

\def\om{\omega}

\def\aut#1{\operatorname{Aut}(#1)}

\def\be{\beta}

\def\al{\alpha}
\def\ga{\gamma}
\def\la{\lambda}
\def\ph{\varphi}
\def\de{\delta}

\def\ggg{\Gamma}

\def\nfs/{NFS}
\def\cdp/{CDP}
\def\cdpz/{CDP${}_0$}
\def\cnL#1{$(L_{#1})$}
\def\cnLl#1{$(L'_{#1})$}

\def\lleq{\preccurlyeq}
\def\ggeq{\succcurlyeq}

\def\cQ{\mathcal{Q}}

\def\cN{\mathcal{N}}
\def\cL{\mathcal{L}}
\def\cp{\mathcal{p}}
\def\cp{\mathcal{S}}
\def\s{{\mathfrak{Q}}}
\def\cG{\mathcal{G}}

\def\gr{{\mathfrak{R}}}

\def\nn{{\mathfrak{N}}}

\def\sq#1#2{(#1)_{#2}}
\def\sqn#1{\sq{#1}{n\in\om}}
\def\sqnn#1{\sqn{#1_n}}

\def\markS#1{$(Q_{#1})$}
\def\llrrd#1#2{\sideset{^{{\rlap{$\scriptstyle{#1}$}\,\,}}}{{#2}}{\operatorname{\D}}}

\def\lrd#1#2{\llrrd{#1}{^{#2}}}
\def\rd#1{\D^{#1}}
\def\ld#1{\llrrd{#1}{}}

\def\lrdd#1#2#3{\llrrd{#1}{^{#2}_{#3}}}
\def\rdd#1#2{\D^{#1}_{#2}}
\def\ldd#1#2{\llrrd{#1}{_{#2}}}
\def\dd#1{\D_{#1}}

\def\lrds#1#2#3{\lrd{#1}{#2}(#3)}

\def\ds#1{\D(#1)}
\def\dps#1#2{\ds{#1;#2}}
\def\dpss#1{\dps{#1}{\s}}

\def\cl#1{\overline{#1}}
\def\clx#1#2{\overline{#2}^{\,#1}}
\def\clv#1{\clx{v}{#1}}
\def\Int{\operatorname{Int}}

\def\id{\operatorname{id}}

\def\bt{\operatorname{\beta}}

\def\oms{\om^*}

\def\limp#1{{\textstyle \lim_{#1}}}

\def\es{\varnothing}
\def\Tau{{\mathcal T}}

\def\nom{{n\in\om}}

\def\sset#1{\{#1\}}
\def\fset#1{\{#1\}}
\def\set#1{\bbset#1\eeset}
\def\bbset#1:#2\eeset{\{#1\,:\,#2\}}

\def\bbsett#1:#2\eesett{\{#1\,:\,\text{#2}\}}

\def\iset#1{\ibbset#1\ieeset}
\def\ibbset#1:#2\ieeset{(#1)_{#2}}

\def\gm{{\mathfrak{M}}}

\def\sod#1{\mathfrak{S}(#1)} 
\def\sse#1{\Expne {#1 \times #1}}

\def\tp{\Tau}
\def\tps{\tp^*}

\def\cP{{\mathcal P}}
\def\cB{{\mathcal B}}

\def\cB{{\mathcal B}}

\def\cU{{\mathcal U}}

\def\cD{{\mathcal D}}

\def\Exp#1{\operatorname{Exp}(#1)}
\def\Expne#1{\operatorname{{Exp}}_*(#1)}

\def\gi{{\mathfrak{i}}}
\def\gm{{\mathfrak{m}}}

\def\til{\tilde}
\def\tD{{\til{\D}}}

\def\wt#1{\widetilde{#1}}

\def\eqdef{\coloneqq}

\def\cX{{\mathcal X}}
\def\cV{{\mathcal V}}
\def\cN{{\mathcal N}}

\def\nov#1{\operatorname{Nov}(#1)}
\def\dev#1{\operatorname{dev}(#1)}
\def\St{\operatorname{St}}
\def\st{\operatorname{st}}

\newcommand\restrA[2]{{
  \left.\kern-\nulldelimiterspace 
  #1 
  \vphantom{\big|} 
  \right|_{#2} 
  }}

\newcommand\restrB[2]{\ensuremath{\left.#1\right|_{#2}}}

\def\restr#1#2{\restrB{#1}{#2}}
\def\sorestr#1#2{\restr{#1}{#2}}

\def\pwr#1_#2{#1^{[#2]}}

\def\ddpsV#1#2{(#1)_{#2}}

\def\D{\Delta}
\def\term#1{{\it #1}}

\def\alp{{\al\in P}}

\def\ppl{{\mathsf{p}}}
\def\pp{{\mathsf{P}}}

\def\dddm#1(#2){N_{#1}(#2)}
\def\dddb#1(#2){B_{#1}(#2)}

\def\wh#1{\widehat{#1}}

\def\et(#1){ (#1)}

\def\gbmt{\ggg^{\wh{\scriptscriptstyle{BM}}}}
\def\godt{\ggg^{\wh{\scriptscriptstyle{OD}}}}
\def\bitem#1,#2.{ $#2\nrightarrow #1$:\ }
\long\def\edemo#1\endedemo{{\rm #1}}

\def\oo#1/{$O_{#1}$}
\def\ooo#1{(\oo{#1}/)}
\def\mgoo#1{O(#1)}
\def\goo#1{$\mgoo{#1}$}

\def\msgoo#1{(#1)}

\def\lrar{\leftrightarrow}

\def\nne{\cN_e}
\def\nnee{\wt\cN_e}
\def\gsod#1{\mathfrak{S}_g(#1)}

\def\rr{$R$}
\def\ep{\varepsilon}
\def\clg#1{\clx{d}{#1}}
\def\clh#1{\clx{h}{#1}}
\def\gppl{{\mathsf{n}}}

\def\gpp{{\mathsf{N}}}
\def\gcP{\cU}
\def\gcQ{\cV}

\def\gD{g\Delta}

\def\gllrrd#1#2{\sideset{^{{\rlap{$\scriptstyle{#1}$}\,\,}}}{{#2}}{\operatorname{\gD}}}

\def\glrd#1#2{\gllrrd{#1}{^{#2}}}
\def\grd#1{\gD^{#1}}
\def\gld#1{\gllrrd{#1}{}}

\def\glrdd#1#2#3{\gllrrd{#1}{^{#2}_{#3}}}
\def\grdd#1#2{\gD^{#1}_{#2}}
\def\gldd#1#2{\gllrrd{#1}{_{#2}}}
\def\gdd#1{\gD_{#1}}

\def\glrds#1#2#3{\glrd{#1}{#2}(#3)}

\def\gds#1{\gD(#1)}
\def\gdps#1#2{\gds{#1;#2}}
\def\gdpss#1{\gdps{#1}{\s}}

\def\grdda#1{\grdd {_{#1}}a}

\def\cnLg#1{$(L^g_{#1})$}

\newtheorem{assertion}{Statement}

\newtheorem{proposition}{Proposition}
\newtheorem{theorem}{Theorem}
\newtheorem*{theorem*}{Theorem}
\newtheorem{lemma}{Lemma}
\newtheorem*{lemma*}{Lemma}
\newtheorem{cor}{Corollary}

\theoremstyle{definition}
\newtheorem{example}{Example}
\newtheorem{definition}{Definition}

\newtheorem{problem}{Problem}
\theoremstyle{remark}
\newtheorem{note}{Remark}
\newtheorem*{note*}{Remark}

\def\setNameGBT#1#2{{\global\expandafter\edef\csname NameGBT{#1}\endcsname{#2}}}
\def\getNameGBT#1{\csname NameGBT{#1}\endcsname}
\def\newlabeldbt#1#2{\setNameGBT{#1}{#2}}
\def\dref#1{\getNameGBT{#1}}
\newlabeldbt{dbmb1}{3}
\newlabeldbt{pbmb4q}{1}
\newlabeldbt{pbmb4}{3}
\newlabeldbt{pbms1}{4}
\newlabeldbt{tbms1}{2}
\newlabeldbt{tbms2}{3}
\newlabeldbt{pgtt1}{5}
\newlabeldbt{pgtt2}{6}
\newlabeldbt{pgtt4}{7}
\newlabeldbt{pgtt3}{8}
\newlabeldbt{psonbd1}{9}
\newlabeldbt{psonbd2}{10}
\newlabeldbt{psonbd3}{11}
\newlabeldbt{psonbd4}{12}
\newlabeldbt{psonbd5}{14}
\newlabeldbt{psonbd6}{15}
\newlabeldbt{psonbd7}{16}
\newlabeldbt{psonbd8}{17}
\newlabeldbt{asonbd1}{1}
\newlabeldbt{psonbd9}{18}
\newlabeldbt{psonbd10}{19}
\newlabeldbt{psonbd11}{20}
\newlabeldbt{pnfs1}{21}
\newlabeldbt{pnfs2}{22}
\newlabeldbt{pnfs3}{23}
\newlabeldbt{pdps1}{24}
\newlabeldbt{adps0}{2}
\newlabeldbt{adps1}{3}
\newlabeldbt{adps2}{4}
\newlabeldbt{pdps2}{25}
\newlabeldbt{pdps3}{26}
\newlabeldbt{adps3}{5}
\newlabeldbt{pdps4}{27}
\newlabeldbt{adps4}{6}
\newlabeldbt{pdps5}{28}
\newlabeldbt{pdps6}{29}
\newlabeldbt{pdps7}{30}
\newlabeldbt{adps5}{7}
\newlabeldbt{adps6}{8}
\newlabeldbt{adps7}{9}
\newlabeldbt{adps8}{10}
\newlabeldbt{pdps9}{31}
\newlabeldbt{pdps10}{32}
\newlabeldbt{pdp1}{33}
\newlabeldbt{adp1}{11}
\newlabeldbt{adp2}{12}
\newlabeldbt{pdp2}{34}
\newlabeldbt{pdp3}{35}
\newlabeldbt{pgd1}{36}
\newlabeldbt{tgd1}{4}
\newlabeldbt{pcdp1}{37}
\newlabeldbt{pcdp2}{38}
\newlabeldbt{pcdp3}{39}
\newlabeldbt{pcdp4}{40}
\newlabeldbt{tcdp1}{5}
\newlabeldbt{pcdp5}{41}
\newlabeldbt{pdsp1}{42}
\newlabeldbt{ptg1}{43}
\newlabeldbt{ptg2}{44}
\newlabeldbt{ptg3}{45}
\newlabeldbt{ctg1}{1}
\newlabeldbt{aqe1}{13}
\newlabeldbt{eqe6}{1}
\newlabeldbt{eqe7}{2}
\newlabeldbt{eqe8}{3}
\newlabeldbt{pqe5}{1}
\newlabeldbt{pqe6}{3}

\begin{document}

\begin{frontmatter}

\title{Continuity in right semitopological groups}
\author{Evgenii Reznichenko} 

\ead{erezn@inbox.ru}

\address{Department of General Topology and Geometry, Mechanics and  Mathematics Faculty, 
M.~V.~Lomonosov Moscow State University, Leninskie Gory 1, Moscow, 199991 Russia}

\begin{abstract}
Groups with a topology that is in consistent one way or another  with the algebraic structure are considered. Classical groups with a topology are topological, paratopological, semitopological, and quasitopological groups. We also study other ways of matching topology and algebraic structure. The minimum requirement in this paper is that the group is a right semitopological group (such groups are often called right topological groups). We study  when a group with a topology is a topological group; research in this direction began with the work of Deane Montgomery and Robert Ellis. (Invariant) semi-neighborhoods of the diagonal are used as a means of study.
\end{abstract}
\end{frontmatter}

\def\tW{\widetilde W}
\def\tG{\widetilde G}
\def\Fp{\Phi}
\def\tc{\Lambda}
\def\tcs{\Lambda^*}

\section{Introduction}


Since the 1936 paper \cite{Montgomery1936} of Montgomery and  the 1957 paper \cite{Ellis1957} of Ellis, the following problem has been studied:
\begin{problem}
Let a topology be given on a group $G$ that is consistent with the group structure of $G$. Under what conditions on the topology and the algebraic structure does the consistency improve? When is $G$ a topological group?
\end{problem}

We study the  following connections between topology and group structure:
\begin{itemize}
\item[\ooo p] multiplication $(g,h)\mapsto gh$ in a group is continuous;
\item[\ooo s] multiplication in a group is separately continuous;
\item[\ooo i] the operation of taking the inverse element $g\mapsto g^{-1}$ is continuous;
\item[\ooo l] left shifts $\la_g : G\to G, x\mapsto gx$ are continuous for all $g\in G$;
\item[\ooo r] right shifts $\rho_g : G\to G, x\mapsto xg$ are continuous for all $g\in G$.
\end{itemize}

Recall that a group $G$ is called
\begin{itemize}
\item
\term{semitopological} if \oo s/  holds;
\item
\term{right semitopological} if \oo r/  holds;
\item
\term{left semitopological} if \oo l/  holds;
\item
\term{paratopological} if \oo p/  holds;
\item
\term{quasitopological} if \oo s/ and \oo i/  hold;
\item
\term{topological} if \oo p/ and \oo i/  holds.
\end{itemize}
In some papers, left (right) semitopological groups are also called left (right) topological  groups.

Montgomery \cite{Montgomery1936} proved that a (locally) complete metrizable semitopological group is a paratopological group and a Polish semitopological group is a topological group. Ellis \cite{Ellis1957} proved that a locally compact Hausdoff paratopological group is a topological group.

Later, the study of the continuity of operations in groups was continued in many papers 
(see \cite{tk2014} and references therein). The most studied objects in this area are semitopological, paratopological and topological groups. The most important class of these is the class of topological groups. The main question in this field is: when is a group  with a topology a topological group?

Groups with a structures weaker than that of semitopological groups  have been studied as well  \cite{Reznichenko1994,SoleckiSrivastava1997,FerriHernandezWu2006,cdp2010,GlasnerMegrelishvili2013,Moors2016} are studied also.

Papers devoted to right topological compact groups  contain also  results on strengthening continuity in such groups. If $G$ is a right topological compact group in which the multiplication $(g,h)\mapsto gh$ is continuous at the identity of $G$, then the inversion map $g\mapsto g^{-1}$ is continuous at the identity group $G$ \cite{ruppert1975}. Theorem \ref{tgmain1} (1) extends this theorem to Baire spaces in the wide class $\cD_d$ defined in Section \ref{sec-gmain}. The class $\cD_d$ includes, for example, locally pseudo-compact spaces and Baire metric spaces.

A right topological group $G$ is called \term{admissible} if there is a dense subset $S$ of $G$ such that $x \mapsto zx$ is continuous for each $z$ in $S$.
We write ``CHART'' for ``compact Hausdorff admissible right topological''.
Let $G$ be a CHART group. Any of the following conditions implies that $G$ is a topological group.
\begin{itemize}
\item[$(C_1)$]
$G$ is metrizable (Theorem 2.1 \cite{namioka1972}).
\item[$(C_2)$]
$G$ is first-countable (Remark after Proposition 1.7 \cite{MoorsNamioka2013}).
\item[$(C_3)$]
$G$ is Fr\'echet (Corollary 8.8 \cite{GlasnerMegrelishvili2013}).
\item[$(C_4)$]
$G$ is tame (Theorem 8.7 \cite{GlasnerMegrelishvili2013}, see also Theorem 3.3 \cite{Moors2016}).
\item[$(C_5)$]
The multiplication of $G$ is feebly continuous. (Proposition 3.2 \cite{Moors2016}, see also Corollary \ref{cgmain1} (2)).
\end{itemize}
Note that the group $G$ in $(C_2)$, $(C_2)$, and $(C_3)$ is metrizable because the compact first-countable and Fr\'echet topological groups are metrizable (Corollary 4.2.2 of \cite{at2009}).
Theorem \ref{tgmain3} (2) strengthens $(C_1)$ and Theorem 3.3 \cite{cdp2010}. Theorem \ref{tgmain1} (2) extends $(C_5)$ to Baire spaces in the class $\cD_d$.
We also obtain the following conditions   for a CHART group to be a topological group:
\begin{itemize}
\item[$(G_6)$]
$G$ has countable $\pi$-character (for example, $G$ is a compact space with countable tightness) (Corollary \ref{cgmain1} (3));
\item[$(G_7)$]
\rm{(MA)}
$w(G)<2^\omega$ (Corollary \ref{cgmain2});
\end{itemize}
Note that condition $(G_6)$ is weaker than $(G_3)$.


To strengthen continuity in groups, the Baire property is needed. There exists a countable metrizable paratopological nontopological group, for example, a countable dense subgroup of the Songenfrey line. In this paper, we study how the Baire-type properties in \cite{rezn2022gbtg} and \cite{rezn2022gbd} affect continuity in topological groups.

In \cite{rezn2022gbd}, Baire-type properties were introduced and studied with the help of semineighborhoods of the diagonal; see Section \ref{sec-tsond}. In right semitopological groups, there are invariant diagonal semineighbourhoods, with the help of which some algebraic-topological properties of right semitopological groups are defined.


\section{Definitions and notation}

The sign $\eqdef$ will be used for equality by definition.

In what follows, we suppose given a group $G$.
Usually, we will denote the identity of $G$ as $e$. For $g\in G$ we denote
\begin{align*}
\la_g&: G\to G,\ x\mapsto gx,
\\
\rho_g&: G\to G,\ x\mapsto xg,
\end{align*}
left and right shifts in $G$.
We denote by $\gm$ and $\gi$ multiplication and inversion in the group:
\begin{align*}
\gm&: G\times G \to G,\ (g,h)\mapsto gh,
\\
\gi&: G\to G,\ g\mapsto g^{-1}.
\end{align*}

Usually, it is assumed that there is some topology on the group $G$. In this case we denote by $\nne$ the family of open neighborhoods of $e$.


The family of all subsets of a set $X$ is denoted by $\Exp X$.
The family of all nonempty subsets of the set $X$ is denoted by $\Expne X$: $\Expne X\eqdef\Exp X \setminus \fset{\es}$.

If $B$ is a subset of the a $A$ then we denote by $B^c=A\setminus B$ the complement to $A$. We use this notation in situations where it is clear from the context which set $A$ is meant.

We define an \term{indexed set} $x=\iset{x_\al:\al\in A}$ as a function on $A$ such that $x(\al)=x_\al$ for $\al\in A$.
If the elements of an indexed set $\cX=\iset{X_\al: \al\in A}$ are themselves sets, then $\cX$ is also called an \term{indexed family of sets}; $\cX$ is a function on $A $: $\cX(\al)=X_\al$ for $\alp$.


We denote by $\aut X$ the set of all homeomorphisms of the space $X$ onto itself.

A subset $M$ of a topological space $X$ is called \term{locally dense} (\term{nearly open} or \term{preopen}) if $M\subset \Int \cl M$.

Let $M\subset X$. If $M$ is the union of a countable number of nowhere dense sets, then $M$ is called a \term{meager} set. Non-meager sets are called \term{sets of the second Baire category} or \term{nonmeager}  sets.
A subset of $M$ is called \term{residual} (\term{comeager}) if $X\setminus M$ is a meager set.

A space $X$ is called  a\term{space of the first Baire category} or a \term{meager space} if the set $X$ is of the first category in the space $X$. A space $X$ is called a \term{space of the second Baire category} or a \term{nonmeager space} if $X$ is not a meager space. A space in which every residual set is dense is called a \term{Baire space}.
A space is nonmeager if and only if some open subspace of it is a Baire space.

A family $\nu$ of nonempty subsets of $X$ is called a \term{$\pi$-net} if for any open nonempty $U\subset X$ there exists $M\in\nu$ such that $M\subset U $.

A $\pi$-network consisting of open sets is called a \term{$\pi$-base}.

A subset $U\subset X$ is called \term{regular open} if $U=\Int{\cl U}$.

A space $X$ is called \term{quasi-regular} if for every nonempty open $U\subset X$ there exists a nonempty open $V\subset X$ such that $\cl V\subset U$.

A space $X$ is called \term{semiregular} if $X$ has a base consisting of regular open sets.

A space $X$ is called \term{$\pi$-semiregular} \cite{Ravsky2001} (or \term{nearly regular} \cite{Ameen2021}) if $X $ is a $\pi$-base consisting of regular open sets.

For a cardinal $\tau$, a set $G\subset X$ is called a \term{set of type $G_\tau$} if $G$ is the intersection of $\tau$ open sets. A space $X$ is called \term{absolute $G_\tau$} if $X$ is of type $G_\tau$ in some compact extension.

A space $X$ is \term{regular at a point} $x\in X$ if for any neighborhood $U$ of the point $x$ there exists a neighborhood $V\ni x$ such that $\cl V\subset U$.

The space $X$ is \term{semiregular at the point} $x\in X$ if there is a base at the point $x$ consisting of regular open sets.

A space $X$ is \term{feebly compact} if any locally finite family of open sets is finite.

For $\ga\subset \Exp X$ and $x\in X$ we denote
\begin{align*}
\St(x,\ga)&\eqdef\set{U\in\ga: x\in \ga},
&
\st(x,\ga)&\eqdef\bigcup\St(x,\ga).
\end{align*}

A space $X$ is called \term{developable} if there exists a sequence of open covers $(\ga_n)_{n\in\om}$ such that for any $x\in X$ the family $\st(x,\ga_n)$ is a base at the point $x$.

A family $\cB$ of open nonempty sets in $X$ is called an \term{outer base of $M\subset X$} if $M\subset U$ for each $U\in\cB$ and for each open $W\supset M$  there exists $U\in \cB$ such that $M\subset U\subset W$.

We denote by $\bt \om$ the space of ultrafilters on $\om$, the Stone-\v{C}ech extension of the discrete space $\om$. We denote by $\oms=\bt\om\setminus \om$  the set of nonprincipal ultrafilters.

Let $\sqnn x$ be a sequence of points in the space $X$  and let $p\in \oms$ be a nonprincipal ultrafilter. A point $x\in X$ is called the \term{$p$-limit} of a sequence $\sqnn x$ if $\set{n\in\om: x_n\in U}\in p$ for any neighborhood $U$  of the point $x$. We will write $x=\limp p x_n=\limp p \sqnn x$ for the $p$-limit $x$.

A space $X$ is said to have \term{countable pseudocharacter} if each point of $X$ is a set of type $G_\de$.

A space $X$ \term{is submetrizable} if there exists a continuous injective mapping of $X$ into a metrizable space.

\section{Weakenings of continuity}\label{sec-vc}

We will use several continuity relaxations following \cite{cdp2010}.
A mapping of topological spaces $f: X\to Y$ is called
\[
\left\{
\text{
\begin{tabular}{c}
continuous
\\
nearly continuous
\\
quasi-continuous
\\
semi-precontinuous
\\
feebly continuous
\end{tabular}
}
\right\}
\text{ at $x\in X$ if }
\left\{
\begin{array}{c}
x\in \Int f^{-1}(V)
\\
x\in \Int \cl{f^{-1}(V)}
\\
x\in \cl{\Int f^{-1}(V)}
\\
x\in \cl{\Int \cl{f^{-1}(V)}}
\\
\Int f^{-1}(V) \neq \es
\end{array}
\right\}
\]
for any neighborhood $V$ of the point $f(x)$. If the properties under consideration are satisfied at every point, then we obtain the definition: the mapping $f$ is called
\[
\left\{
\text{
\begin{tabular}{c}
continuous
\\
nearly continuous
\\
quasi-continuous
\\
semi-precontinuous
\\
feebly continuous
\end{tabular}
}
\right\}
\text{ if }
\left\{
\begin{array}{c}
f^{-1}(V)\subset \Int f^{-1}(V)
\\
f^{-1}(V)\subset \Int \cl{f^{-1}(V)}
\\
f^{-1}(V)\subset \cl{\Int f^{-1}(V)}
\\
f^{-1}(V)\subset \cl{\Int \cl{f^{-1}(V)}}
\\
\Int f^{-1}(V) \neq \es
\end{array}
\right\}
\]
for any open nonempty $V\subset X$.

We call a mapping $f$ \term{a feebly homeomorphism} if $f$ is a bijection and the mappings $f$ and $f^{-1}$ are feebly continuous.

We call a mapping $f$ \term{a quasi homeomorphism} if $f$ is a bijection and the mappings $f$ and $f^{-1}$ are quasi continuous.

Let $X$, $Y$, $Z$ be three topological spaces, and let $f : X \times Y \to Z$ be a function
from $X \times Y$ into $Z$. Recall that $f$ is called \term{quasi-continuous with respect to the second
variable} \cite{Piotrowski1980}
at $(x, y)$ if for every open neighborhood $W$ of $f (x, y)$ and every open neighborhood
$U \times V$ of $(x, y)$, there are an open neighborhood $V'$ of $y$ and a nonempty open set
$U' \subset U$ such that $f (U' \times V') \subset W$. Quasi-continuity with respect to the second
variable plays an important role in the theory of separate vs. joint continuity. It
is called \term{strong quasi-continuity} in \cite{Bouziad1996} and \cite{KenderovKortezovMoors2001}, where it is applied to the study
of the problem when a semitopological group is a topological group.

Similarly, define that $f$ is called \term{quasi-continuous with respect to the first variable} at $(x, y)$ if for every open neighborhood $W$ of $f (x, y)$ and every open neighborhood
$U \times V$ of $(x, y)$, there are an open neighborhood $U'$ of $x$ and a nonempty open set
$V' \subset V$ such that $f (U' \times V') \subset W$.

Let us list the connections between the topology and the group structure of the group $G$.
\begin{itemize}
\item[\ooo {t}] holds \ooo p and \ooo i, that is, $G$ is a topological group;
\item[\ooo {pe}] multiplication $\gm$ in a group is continuous in $(e,e)$;
\item[\ooo {qpe}] $\gm$ is quasi-continuous in the first coordinate in $(e,e)$;
\item[\ooo {sqpe}] $\gm$ is quasi-continuous in $(e,e)$;
\item[\ooo {fpe}] $\gm$ is feebly continuous in $(e,e)$;
\item[\ooo {fi}] the operation of taking the inverse element $\gi$ is feebly continuous;
\item[\ooo {ie}] $\gi$ is continuous at the identity $e$;
\item[\ooo {nie}] $\gi$ is almost continuous in $e$;
\item[\ooo {qie}] $\gi$ is quasi-continuous in $e$;
\item[\ooo {sie}] $\gi$ is semi-continuous in $e$;
\item[\ooo {fie}] $\gi$ is feebly continuous in $e$;
\item[\ooo {fl}] left shifts of $\la_g$ are feebly continuous for any $g\in G$;
\item[\ooo {dfl}] there exists a dense $H\subset G$ such that $\la_g$ are feebly continuous for any $g\in H$;
\item[\ooo {dfl^*}] there exists a dense $H\subset G$ such that $\la_g$ are feebly homeomorphisms for any $g\in H$;
\item[\ooo {ql}] left shifts of $\la_g$ are quasi continuous for any $g\in G$;
\item[\ooo {dql}] there exists a dense $H\subset G$ such that $\la_g$ are quasi continuous for any $g\in H$;
\item[\ooo {dql^*}] there exists a dense $H\subset G$ such that $\la_g$ are quasi homeomorphisms for any $g\in H$;
\item[\ooo {dl}] there is a dense $H\subset G$ such that $\la_g$ are continuous for any $g\in H$;
\item[\ooo {dl^*}] there is a dense $H\subset G$ such that $\la_g$ are homeomorphisms for any $g\in H$;
\end{itemize}

We set
\begin{align*}
\tc_{f}(G)&\eqdef\set{g\in G:\text{$\la_g$ is feebly continuous}},
& \tcs_{f}(G)&\eqdef\set{g\in G:g,g^{-1}\in \tc_{f}(G)},
\\
\tc_{q}(G)&\eqdef\set{g\in G:\text{$\la_g$ is quasi continuous}},
& \tcs_{q}(G)&\eqdef\set{g\in G:g,g^{-1}\in \tc_{q}(G)},
\\
\tc(G)&\eqdef\set{g\in G:\text{$\la_g$ is continuous}},
& \tcs(G)&\eqdef\set{g\in G:g,g^{-1}\in \tc(G)}.
\end{align*}

Easy to check
\begin{assertion}\label{avc1}
The sets $\tcs_{f}(G)$, $\tcs_{q}(G)$, $\tcs(G)$ are subgroups of the group $G$.
The following alternative definitions of the listed properties are true:
\ooo {l} $\tc(G)=G$;
\ooo {fl} $\tc_f(G)=G$;
\ooo {ql} $\tc_q(G)=G$;
\ooo {dl} $\cl{\tc(G)}=G$;
\ooo {dfl} $\cl{\tc_f(G)}=G$;
\ooo {dql} $\cl{\tc_q(G)}=G$;
\ooo {dl^*} $\cl{\tcs(G)}=G$;
\ooo {dfl^*} $\cl{\tcs_f(G)}=G$;
\ooo {dql^*} $\cl{\tcs_q(G)}=G$.
\end{assertion}

Let us give an alternative definition of some of the listed properties. Let $\tps$ be all nonempty open subsets of $G$.
\begin{itemize}
\item[\ooo {pe}] for any $U\in\nne$ there exists $V\in\nne$ so $V^2\subset U$;
\item[\ooo {sqpe}] for any $U\in\nne$ there are $V\in\nne$ and $W\in\tps$ so $VW\subset U$;
\item[\ooo {fpe}] for any $U\in\nne$ there are $V\in\tps$ and $W\in\tps$ so $VW\subset U$;
\item[\ooo {fi}] $\Int U^{-1}\neq \es$ for any $U\in\tps$;
\item[\ooo {ie}] $e\in\Int U^{-1}$ for any $U\in\nne$;

\item[\ooo {nie}] $e\in \Int \cl{U^{-1}}$ for any $U\in\nne$;
\item[\ooo {qie}] $e\in \cl{\Int U^{-1}}$ for any $U\in\nne$;
\item[\ooo {sie}] $e\in \cl{\Int \cl{U^{-1}}}$ for any $U\in\nne$;

\item[\ooo {fie}] $\Int U^{-1}\neq \es$ for any $U\in\nne$;
\item[\ooo {l}] $gU\in\nne$ for any $U\in\tps$ and $g\in G$;
\item[\ooo {r}] $Ug\in\nne$ for any $U\in\tps$ and $g\in G$;
\item[\ooo {fl}] $\Int gU\neq\es$ for any $U\in\tps$ and $g\in G$;
\item[\ooo {ql}] $g\in\cl{\Int Ug}$ for any $U\in\tps$ and $g\in G$.
\end{itemize}

\begin{definition}\label{dvc1}
For
\begin{align*}
&\sset{x_1,x_2,...,x_n} \subset
\\
&\sset{p,s,i,l,r,pe,qpe,sqpe,fpe,ie,nie,fie,qie,sie,fl,dfl, ql,dql,dl, dfl^*,dql^*,dl^*}
\end{align*}
we say that a group $G$ is \goo{x_1,x_2,...,x_n}-topological if the condition \ooo{x_i} is satisfied for all $i=1,2,...,n$.
\end{definition}

Topological groups are \goo{p,i}-topological groups, and other introduced classes of groups with topology can be defined similarly.


\begin{proposition}\label{pvc1}
Below, $\msgoo A\to \msgoo B$ means that \goo A-topological group is \goo B-topological group
and $\msgoo A\lrar \msgoo B$ means that $G$ is a \goo A-topological group if and only if $G$ is a \goo B-topological group.
{\def\s{\qquad}
\begin{gather*}
\msgoo{p,i}\lrar\msgoo{t}
\s
\msgoo{l,r}\lrar\msgoo{s}
\s
\msgoo{s,pe}\lrar\msgoo{p}
\s
\msgoo{s,ie}\lrar\msgoo{s,i}
\\
\msgoo{l,i}\lrar\msgoo{r,i}\lrar\msgoo{s,i}
\\
\msgoo{p}\to\msgoo{s}
\s
\msgoo{p}\to\msgoo{pe}\to\msgoo{sqpe}
\s
\msgoo{i}\to\msgoo{ie}\to\msgoo{fie}
\\
\msgoo{l}\to\msgoo{fl}\to\msgoo{dfl}
\\
\msgoo{ie}\to\msgoo{nie}\to\msgoo{sie}
\s
\msgoo{ie}\to\msgoo{qie}\to\msgoo{sie}
\s
\msgoo{qie}\to\msgoo{fie}
\\
\msgoo{dfl^*,fie,r}\to\msgoo{fl,fi}
\s
\msgoo{dfl^*,fpe,r}\to\msgoo{sqpe}
\\
\msgoo{p,sie}\to\msgoo{t}
\s
\msgoo{dfl^*,ie,pe,r}\to\msgoo{t}
\end{gather*}}
\end{proposition}
\begin{proof}
Most of these statements are either trivial or widely known (see \cite{cdp2010} and \cite{at2009}). Let us prove new nontrivial assertions.

$\msgoo{p,sie}\to\msgoo{t}$. Theorem 2.3 from \cite{cdp2010}, Lemma 1.2 from \cite{arh-rezn2005}.

$\msgoo{dfl^*,fie,r}\to\msgoo{fl,fi}$.
Let us show that $\msgoo{dfl^*,fie,r}\to\msgoo{fi}$. Let $U\in\tps$.
Take $g\in \tcs_f(G)\cap U$. From
\[
U^{-1}=g^{-1}(Ug^{-1})^{-1}
\]
it follows that $\Int U^{-1}\neq\es$. 
Let us show $\msgoo{fl}$. We have $hU=(U^{-1} h^{-1})^{-1}$ for any $h\in G$. Since $\msgoo{fi,r}$, then $\Int hU\neq\es$.

$\msgoo{dfl^*,ie,pe,r}\to\msgoo{t}$. It follows from the previous subsection that $G$ is a \goo{fl}-topological group.
By virtue of $\msgoo{s,pe}\lrar\msgoo{p}$ and $\msgoo{s,ie}\lrar\msgoo{s,i}$,
it suffices to prove that $G$ is a \goo{l}-topological group. Let $g\in G$ and $U\in\tps$. You need to check that $gU$ is open. Let $q\in U$. It suffices to check that $gq\in\Int gU$, which is equivalent to $g\in \Int gV$, where $V=Uq^{-1}\in \nne$.
There is $W\in \nne$ for which $WW^{-1}\subset V$. Let $s\in \Int gW$ and $w=g^{-1}s\in W$. There is $S\in\nne$ so $Ss\subset gW$. Then $s=gw$, $g=sw^{-1}$ and
\[
Sg=Ssw^{-1}\subset gWw^{-1}\subset gWW^{-1}\subset gV.
\]
Hence $g\in \Int gV$.

$\msgoo{dfl^*,fpe,r}\to\msgoo{sqpe}$. Let $U\in \nne$. Then $WVx\subset U$ for some $x\in G$, $W\in\tps$ and $V\in\nne$. Let $y\in \tcs_f(G) \cap W$. Then $Sy\subset W$ for some $S\in\nne$.
For some $z\in \Int yV$ and $Q\in\nne$, $Qz\subset yV$. We get $SQzx \subset U$.
\end{proof}

\section{\rr-topological groups}\label{sec-rtg}

Let us call the group $G$
\begin{itemize}
\item
\term{\rr-semitopological} if $G$ is a \goo{r}-topological group, that is, $G$ is a right semitopological group;
\item
\term{\rr-paratopological} if $G$ is a \goo{pe,r}-topological group;
\item
\term{\rr-quasitopological} if $G$ is a \goo{ie,r}-topological group;
\item
\term{\rr-topological} if $G$ is a \goo{ie,pe,r}-topological group.
\end{itemize}

Let $N$ be a real-valued function on $G$. We shall call $N$ a \term{prenorm} (\cite{at2009}, Section 3.3) on
$G$ if the following conditions are satisfied for all $x, y \in G$:
\begin{itemize}
\item[(PN1)]
$N(e) = 0$;
\item[(PN2)]
$N(xy) \leq N(x) + N(y)$;
\item[(PN3)]
$N(x^{-1}) = N(x)$.
\end{itemize}
If, in addition, the condition
\begin{itemize}
\item[(PN4)]
$N(x) \neq 0$ for $x\neq e$,
\end{itemize}
then $N$ is called \term{norm}.
In \cite{BinOst2008} the prenorm is called the pseudonorm.

\begin{assertion}[Proposition 3.3.1 and Proposition 3.3.2 from \cite{at2009}]\label{artg1}
Let $N$ be a prenorm on the group $G$. Then for $x,y\in G$
\begin{itemize}
\item
$N(x)\geq 0$;
\item
$|N(x)-N(y)|\leq N(xy^{-1})$.
\end{itemize}
\end{assertion}

\begin{assertion}\label{artg2}
Let $\sqnn U$ be a sequence of subsets of the group $G$such that $U_n=U_n^{-1}$ and $U_{n+1}^2\subset U_n$ for $\nom$.
\begin{itemize}
\item[{\rm (1)}]
There is a prenorm $N$ on $G$ such that
\[
\set{x\in G: N(x)< 1/2^n} \subset U_n \subset \set{x\in G: N(x)\leq 2/2^n}
\]
for $\nom$.
\item[{\rm (2)}] If $G$ is an \rr-semitopological group and $e\in\Int U_n$ for $\nom$, then the prenorm $N$ is continuous.
\end{itemize}
\end{assertion}
\begin{proof}
Condition (1) is actually proved in Lemma 3.3.10 from \cite{at2009}. Let us prove (2). Let $y\in G$ and $\ep>0$. There is $\nom$ for which $2/2^n<\ep$. Then $y\in \Int U_n y$. Let $x\in U_n y$. Then $xy^{-1}\in U_n$ and
\[
|N(x)-N(y)|\leq N(xy^{-1}) \leq 2/2^n < \ep.
\]
\end{proof}
For the prenorm $N$ we denote
\[
B_N\eqdef \set{g\in G: N(g)<1}.
\]

\begin{proposition}\label{prtg1}
Let $G$ be an \rr-topological group and $U\in\nne$. Then $e\in B_N\subset U$ for some continuous prenorm $N$ on $G$.
\end{proposition}
\begin{proof} There exists $\sqnn U$ there is a sequence of subsets of the group $G$such that
$U=V_0$, $V_n\in\nne$ and $V_{n+1}^2\subset V_n$ for $\nom$. Let $U_n=V_n\cap V_n^{-1}$. The prenorm $N$ from Statement \ref{artg2} is the desired one.
\end{proof}

Denote
\[
d_N(x,y) \eqdef N(xy^{-1})\text{ and }B_N(x,r)\eqdef \sset{y\in G: d_N(x,y)<r}
\]
for $x,y\in G$ and $r>0$. The function $d_N$ is a right-invariant peseudometric on $G$ and $d_N$ is continuous if the prenorm $N$ is continuous. The right-handed metric $d$ defines the prenorm $N_d$: $N_d(x)=d(e,x)$. Note that $d=d_{N_d}$. Thus, there is a natura one-to-one correspondence between prenorms and right-invariant metrics on $G$, and continuous right-invariant pseudometrics correspond to continuous prenorms.

A family $\set{d_\al:\al\in A}$ of pseudometrics on a space $X$ defines the topology of the space $X$ if the family of sets open with respect to $d_\al$ for $\al\in A$ form a prebase $X $.

\begin{proposition}\label{prtg2}
Let $G$ be a topological group. A group $G$ is an \rr-topological group if and only if the topology of $G$ is given by a family of right-invariant pseudometrics.
\end{proposition}
\begin{proof}
Let $G$ be an \rr-topological group.
It follows from Proposition \ref{prtg1} that $U\in \nne$ has a continuous prenorm $N_U$ for which $B_{N_U}\subset U$. Let $d_U=d_{N_U}$. The family $\set{d_U: U\in\nne}$ of continuous right-invariant metrics defines the topology $G$.

Let $\set{N_\al:\al\in A}$ be a family of continuous prenorms such that
$\set{d_{N_\al}:\al\in A}$ is a family of right-invariant metrics defining the $G$ topology. Right shifts $\rho_g$ are continuous in $G$.
Let $U\in\nne$. There is $\ep>0$ and a finite $B\subset A$ such that
$\bigcap_{\al\in B} B_{N_\al}(e,\ep)\subset U$.
Let $V=\bigcap_{\al\in B} B_{N_\al}(e,\ep/2)$. Then $V\in \nne$, $V=V^{-1}$ and $V^2\subset U$.
\end{proof}


\begin{proposition}\label{prtg3}
Let $G$ be a $T_0$ \rr-topological group.
\begin{itemize}
\item[{\rm (1)}] If $G$ is a first countable space, then there exists a continuous right-invariant metric $d$ on $G$ that defines a topology on $G$.
\item[{\rm (2)}] If $G$ is a space of countable pseudocharacter, then there exists a continuous right-invariant metric $d$.
\end{itemize}
\end{proposition}
\begin{proof} Let $\sqnn W$ be a sequence of neighborhoods of unity, which in case (1) forms a base in $e$ and $\sset e = \bigcap\nom W_n$ in case (2). There is a sequence $\sqnn V$ of neighborhoods of unitysuch that $V_{n+1}^2 \subset V_n\subset W_n$ for $\nom$. Let $U_n=U_n\cap U_n^{-1}$ for $\nom$. Let $N$ be the prenorm as in Statement \ref{artg2} and $d=d_N$.
\end{proof}

\begin{theorem}\label{trtg1}
Let $G$ be a $T_0$ \rr-topological group.
\begin{itemize}
\item[{\rm (1)}] A group $G$ is metrizable if and only if $G$ is first countable.
\item[{\rm (2)}] A group $G$ is submetrizable if and only if $G$ has a countable pseudocharacter.
\end{itemize}
\end{theorem}

Theorem \ref{trtg1} is a generalization of the Birkhoff-Kakutani theorem, see Theorem 3.3.12 and Theorem 3.3.16 from \cite{at2009}.

Proposition \ref{pvc1} implies
\begin{theorem}\label{trtg2}
Let $G$ be an \rr-topological group. If the set $\tcs_f(G)$ dense in $G$, then $G$ is a topological group.
\end{theorem}

\begin{example}\label{ertg1}
Let us describe Example (d) from \cite{milnes1993}. Let $\T=\set{z\in \C: |z|=1}$, $\ph: \T\to \T$ be a discontinuous automorphism of the group $\T$ such that $\ph^2=\ id_\T$. We set $H=\sset{\id_\T,\ph}$ with discrete topology, $G=\T\rtimes H$, $G$ as a set, and the space is homeomorphic to $\T\times H$. Recall the definition of multiplication in a semidirect product $G$:
\[
(t_1,h_1)\cdot (t_2,h_2)=(t_1 h_1(t_2), h_1h_2).
\]
The group $G$ is a compact metrizable right semitopological group. The subgroup $\T$ is normal and clopen in $G$. Since $\T$ is a topological group, $G$ is a \goo{ie,pe}-topological group. Finally, $G$ is a metrizable compact \rr-topological group which is not a topological group.
\end{example}

\section{Properties of spaces defined by semi-neighborhoods of the diagonal}
\label{sec-tsond}
\subsection{Semi-neighborhoods of the diagonal}\label{sec-sonbd}

Let $X$ be a set, $P,Q\subset X\times X$. Denote
\begin{align*}
P(x)&\eqdef\set{y\in X: (x,y)\in P},
\\
P(M)&\eqdef\set{y\in X: \text{ exists }x\in M \text{ so }(x,y)\in P}=\bigcup_{x\in M}P (x),
\\
\sorestr PM &\eqdef P\cap (M\times M),
\\
P^{-1}&\eqdef\set{(x,y):(y,x)\in P},
\\
P\circ Q&\eqdef\set{(x,y):\text{there is }z\in X\text{ so }(x,z)\in P\text{ and }(z,y)\in Q}
\end{align*}
for $x\in X$ and $M\subset X$.

\begin{definition}
Let $X$ be a space, $P\subset X\times X$. We call a set $P$ \term{semiopen} if every $P(x)$ is open.
We call $P$ a \term{diagonal semi-neighbourhood} if $x \in \Int P(x)$ for all $x\in X$.
We denote by $\sod X$ the family of all semi-neighborhoods of the diagonal of the space $X$.
\end{definition}

\begin{note}
In \cite{rezn2008}, in the definition of a semi-neighbourhood and a diagonal, it was assumed that the set $P(x)$ is open, that is, $P$ was assumed to be semi-open. The concept of a semi-open semi-neighborhood of the diagonal essentially coincides with the concept of \cite{mtw2007} neighbourhood assignment: a map $P:X\to \tp$ is called a \term{neighbourhood assignment} if $x\in P(x)$ and $ P(x)$ is open to $x\in X$.

A set $P$ is a semi-neighbourhood of the diagonal if and only if $P$ contains some semi-open semi-neighborhood of the diagonal.
\end{note}

If $X$ is a topological space, then we denote
\begin{align*}
\clv P&\eqdef \set{(x,y)\in X\times X: y\in \cl{P(x)}}.
\end{align*}
It is clear that $\clv P(x)=\cl{P(x)}$ for $x\in X$.

On $\Exp{\sod X}$ we introduce an order relation and an equivalence relation related to it. For $\cP,\cQ\subset \sod X$ we set
\begin{itemize}
\item
$\cP\lleq \cQ$ if and only if for any $Q\in \cQ$ there exists $P\in\cP$such that $P\subset Q$;
\item
$\cP\sim \cQ$ if and only if $\cP\lleq \cQ$ and $\cQ\lleq \cP$.
\end{itemize}
Note that if $\cP\supset \cQ$, then $\cP\lleq \cQ$.
For $\cP\subset \sod X$ we denote
\[
\pp^e(\cP)\eqdef\set{R\in \sod X: P\subset R,\text{ for some }P\in\cp}.
\]

For $\cP\subset \sod X$ we denote
\begin{align*}
\pp(\cP)&\eqdef\cP,
\\
\pp_+(\cP)&\eqdef\set{Q\circ P: Q\in\sod X\text{ and }P\in \cP},
\\
\pp^v(\cP)&\eqdef\set{\clv P: P\in\cP},
\\
\pp^c(\cP)&\eqdef\set{\cl P: P\in\cP},
\\
\pp_+^c(\cP)&\eqdef \pp^v(\pp_+(\cP)).
\end{align*}

\begin{proposition}[Proposition \dref{psonbd3} \cite{rezn2022gbd}] \label{psonbd3}
Let $X$ be a space.
For $\cP \subset \sod X$
\begin{itemize}
\item
$\cP\lleq \pp_+(\cP)\lleq \pp_+^c(\cP)$
\item
$\cP\lleq \pp^v(\cP)\lleq \pp^c(\cP)\lleq \pp_+^c(\cP)$
\end{itemize}
\end{proposition}

Let $\la$ be an ordinal, $P_\al\subset X\times X$ for $\al<\la$. By induction on $\be<\ga$ we define the sets $\ppl_\be(\cp_\be),\ppl^c_\be(\cp_\be)\subset X\times X$, where $\cp_\be= (P_\al)_{\al<\be}$.

$\be=0$. We set
\[
\ppl_0(\cp_0)\eqdef\ppl_0^c(\cp_0)\eqdef\sset\es.
\]

$\be=1$. We set
\[
\ppl_1(\cp_1)\eqdef P_0,
\qquad
\ppl_1^c(\cp_1)\eqdef \cl{P_0}.
\]

$1<\be\leq \ga$. We set
\begin{align*}
\ppl_\be(\cp_\be)&\eqdef\begin{cases}
Q_\be,&\text{if $\be$ is cardinal limit}
\\
Q_\be\circ P_{\be'},&\text{if $\be=\be'+1$}
\end{cases}
\\
\ppl^c_\be(\cp_\be)&\eqdef\begin{cases}
\cl{Q^c_\be},&\text{if $\be$ is cardinal limit}
\\
\cl{Q^c_\be\circ P_{\be'}},&\text{if $\be=\be'+1$}
\end{cases}
\end{align*}
where
\[
Q_\be=\bigcup_{\al<\be} \ppl_\al(\cp_\al)
\qquad
\text{and}
\qquad
Q^c_\be=\bigcup_{\al<\be} \ppl^c_\al(\cp_\al).
\]

For $\cP\subset \sod X$ we denote
\begin{align*}
\pp_\ga(\cP)&\eqdef\set{\ppl_\ga(\cp): \cp\in \cP^\ga},
\\
\pp^c_\ga(\cP)&\eqdef\set{\ppl^c_\ga(\cp): \cp\in \cP^\ga},
\\
\sorestr \cP Y&\eqdef\set{\sorestr PY: P\in\cP}
\end{align*}
for $Y\subset X$.


\begin{proposition}[Proposition \dref{psonbd4} \cite{rezn2022gbd}] \label{psonbd4}
Let $X$ be a space, $\cP \subset \sod X$, and $1<\al<\be$ ordinals.
Then
\begin{align*}
\pp_+(\cP) \lleq \pp_2(\cP) \lleq &\pp_\al(\cP) \lleq \pp_\be(\cP) \lleq \pp_\be^c(\cP)
\\
\pp^c(\cP) = \pp_1^c(\cP) \lleq &\pp_\al^c(\cP) \lleq \pp_\be^c(\cP)
\\
\pp^c(\cP) \lleq \pp^c_+(\cP) \lleq &\pp_2^c(\cP)
\end{align*}
\end{proposition}

\subsection{Normal square functors}\label{sec-nfs}

A functor $\s$ that associates a topological space $X$ with a family $\s(X)\subset \sse X$ of nonempty subsets of $X\times X$ is called the \term{square functor}. A square functor $\s$ is called a \term{normal square functor} (\nfs/) if the following conditions are satisfied:
\begin{itemize}
\item[\markS1] if $f:X\to Y$ is a homeomorphism, then
\[
\s(Y)=\set{(f\times f)(P): P\in \s(X)};
\]
\item[\markS2] if $U$ is an open nonempty subset of $X$, then
\[
\s(U)=\sorestr {\s(X)}U;
\]
\item[\markS3] if $S\in \s(X)$, $S\subset Q \subset X\times X$, then $Q\in \s(X)$;
\item[\markS4] if $S\in \s(X)$ then $\cl S = X\times X$.
\end{itemize}

We introduce an order relation on normal square functors. For \nfs/ $\s$ and $\gr$ $\s\lleq\gr$ is true if and only if $\s(X)\supset \gr(X)$ for any space $X$.

Let us define \nfs/ which we will use. Let $k\in\sset{d,h,v,s,a}$, $X$ be a space, $S\in\sse X$. Then $S\in \s_k(X)$ if and only if the condition \cnL k:
\begin{itemize}
\item[\cnL d] $\cl S=X\times X$;
\item[\cnL v] $\cl{S(x)}=X$ for any $x\in X$;
\item[\cnL h] $\cl{S^{-1}(x)}=X$ for any $x\in X$;
\item[\cnL s] $M\times X\subset S$ for some $M\subset X=\cl M$;
\item[\cnL a] $S=X$.
\end{itemize}


\begin{proposition}[Proposition \dref{pnfs1} \cite{rezn2022gbd}] \label{pnfs1}
Let $k\in\sset{d,h,v,s,a}$, $X$ be a space and $S\in \sse X$. Then $S\in\s_k(X)$ if and only if for any $x\in X$ and any neighborhood $U\subset X$ of $x$ the condition \cnLl k is satisfied:
\begin{itemize}
\item[\cnLl d] $\cl{S(U)}=X$;
\item[\cnLl v] $\cl{S(x)}=X$;
\item[\cnLl h] $S(U)=X$;
\item[\cnLl s] $S(z)=X$ for some $z\in U$;
\item[\cnLl a] $S(x)=X$.
\end{itemize}
\end{proposition}

\begin{proposition}[Proposition \dref{pnfs2} \cite{rezn2022gbd}] \label{pnfs2}
For any \nfs/ $\s$
\begin{align*}
&\s_d\lleq\s\lleq\s_a,
\\
&\s_h\lleq\s_s.
\end{align*}
\end{proposition}

We introduce an increment operation on \nfs/: for \nfs/ $\s$ we define \nfs/ $\s^+$. Let $X$ be a space, $Q\in\sse X$. Then $Q\in\s^+(X)$ if and only if $P\circ S\subset Q$ for some $P\in\sod X$ and $S\in\s(X)$.

\begin{proposition}[Proposition \dref{pnfs3} \cite{rezn2022gbd}] \label{pnfs3}
\begin{align*}
& \s_v \lleq \s^+_d
\\
&\s_a = \s^+_s = \s^+_h
\end{align*}
\end{proposition}

\subsection{\texorpdfstring{$\dps \cP\s$}{DPS}-Baire spaces}\label{sec-dps}

\begin{definition}
Let $\s$ be a normal square functor, $(X,\tp)$ a space,
$\tps=\tp\setminus\sset\es$, $\cP\subset \sod X$. Let us call $X$
\begin{itemize}
\item \term{$\dpss \cP$-nonmeager} ($\dpss \cP$-nonmeager) space if for any $P\in\cP$ there exists $V\in\tps$, so $\sorestr P V\in \s(V)$;
\item \term{$\dpss \cP$-Baire space} ($\dpss \cP$-Baire space) if for any $P\in\cP$ and any $U\in\tps$ there exists an open nonempty $ V\subset U$, so the condition $\sorestr P V\in \s(V)$ is satisfied.
\end{itemize}
\end{definition}

\begin{proposition}[Proposition \dref{pdps7} \cite{rezn2022gbd}] \label{pdps7}
Let $X$ be a space, $\cP \subset \sod X$.
In the diagrams below, the arrow
\[
\ddpsV Fk \to \ddpsV Gl
\]
means that $F,G:\Exp{\sod X}\to \Exp{\sod X}$ are mappings $k,l\in\sset{d,v,h,s,a}$ and terms
\begin{itemize}
\item[\rm (1)] If $X$ is a $\dps{F(\cP)}{\s_k}$-nonmeager space, then $X$ is $\dps{G(\cP)}{\s_l}$-nonmeager space and
\item[\rm (2)] if $X$ is a $\dps{F(\cP)}{\s_k}$-Baire space, then $X$ is $\dps{G(\cP)}{\s_l}$-Baire space.
\end{itemize}

{
\def\xp#1{\ddpsV{\pp}#1}
\def\xpu#1#2{\ddpsV{\pp^{#1}}#2}
\def\xpl#1#2{\ddpsV{\pp_{#1}}#2}
\def\x#1#2#3{\ddpsV{\pp_{#1}^{#2}}#3}
\def\xc#1#2{\ddpsV{\pp_{#1}^{c}}#2}
\def\rr{\arrow[r,leftrightarrow]}

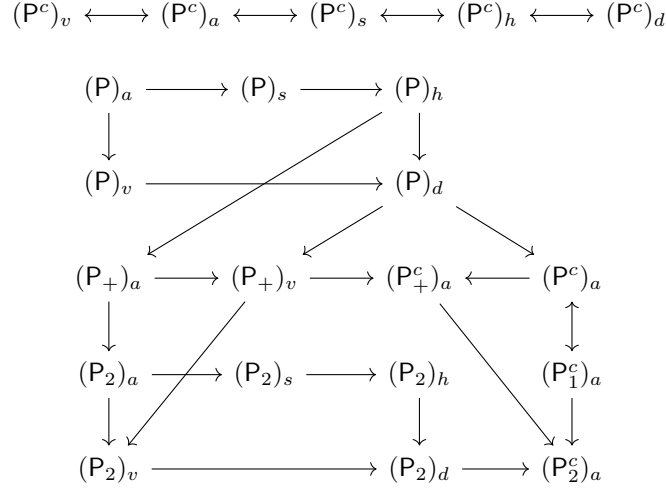
\[
\begin{tikzcd}
\xpu cv \rr&\xpu ca \rr&\xpu cs \rr&\xpu ch \rr&\xpu cd
\end{tikzcd}
\]

\[
\begin{tikzcd}
\xp a \arrow[r]\arrow[d] & \xp s \arrow[r] & \xp h \arrow[d] \arrow[ddll]& \\
\xp v \arrow[rr] && \xp d \arrow[dl]\arrow[dr]&\\
\xpl+a \arrow[r]\arrow[d] & \xpl+v \arrow[r] \arrow[ddl] & \x+ca \arrow[ddr] &\xpu ca \arrow[l] \arrow [d,leftrightarrow]\\
\xpl 2a \arrow[r]\arrow[d] & \xpl 2s \arrow[r] & \xpl 2h \arrow[d] & \xc 1a \arrow[d]\\
\xpl 2 v \arrow[rr] && \xpl 2 d \arrow[r] & \xc 2a\\
\end{tikzcd}
\]

Let $1<\al<\be$ be ordinals.

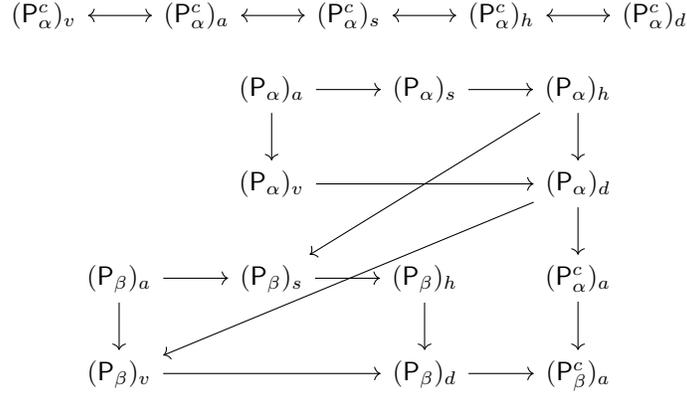
\def\xca#1{\xc{\al}{#1}}
\[
\begin{tikzcd}
\xca v \rr&\xca a \rr&\xca s \rr&\xca h \rr&\xca d
\end{tikzcd}
\]

\[
\begin{tikzcd}
&\xpl \al a \arrow[r]\arrow[d] & \xpl \al s \arrow[r] & \xpl \al h \arrow[ddll] \arrow[d] \\
&\xpl \al v \arrow[rr] && \xpl \al d \arrow[ddlll] \arrow[d]\\
\xpl \be a \arrow[r]\arrow[d] & \xpl \be s \arrow[r] & \xpl \be h \arrow[d] & \xc \al a \arrow[d]\\
\xpl \be v \arrow[rr] && \xpl \be d \arrow[r] & \xc \be a\\
\end{tikzcd}
\]
}

\end{proposition}

\subsection{\texorpdfstring{$\ds \s$}{DS}-Baire spaces}\label{sec-dp}

\begin{definition}
Let $\s$ be a normal square functor, $X$ be a space, $\ga$ be an ordinal,
$\cP = \sod X$. Let $\pp^k_l: \Exp{\sod X}\to \Exp{\sod X}$ be one of the mappings considered in the \ref{sec-sonbd} section:
\[
\pp^k_l\in\sset{\pp,\pp^e,\pp^v,\pp^c,\pp_+,\pp_+^c,\pp_\ga,\pp_\ga^c}.
\]
We say that $X$ is a \term{$\lrds kl{\s}$-Baire} (\term{$\lrds kl{\s}$-nonmeager}) space if $X$ is $\dps{ \pp^k_l(\sod X)}\s$-Baire ($\dps{\pp^k_l(\sod X)}\s$-nonmeager) space. For
\[
\tD\in\sset{\D,\ld e,\ld v,\ld c,\rd +,\lrd c+,\rd\ga,\lrd c\ga}
\]
we have defined $\tD(\s)$-Baire ($\tD(\s)$-nonmeager) spaces. For
\[
k\in\sset{d,h,v,s,a}
\]
we say $X$ is a \term{${\tD}_k$-Baire} (\term{${\tD}_k$-nonmeager}) space if $X$ is $\tD(\s_k)$ -Baire ($\tD(\s_k)$-nonmeager) space. For
\[
\til{\tD}\in\sset{\dd k, \ldd {e}k,\ldd {v}k,\ldd {c}k,\rdd {+}k,\lrdd c{_+} k,\rdd {_\ga}k,\lrdd c{_\ga}k}
\]
we have defined $\til{\tD}$-Baire ($\til{\tD}$-nonmeager) spaces. Also, if the subscript is not written, $d$ is implied: $\lrd kl$-nonmeager ($\lrd kl$-Baire) is $\lrdd kld$-nonmeager ($\lrdd kld$-Baire). A separate direct definition for the most important classes of spaces: the space $X$ is called
\begin{itemize}
\item[\rm (1)] $\D$-nonmeager ($\D$-Baire) if $X$ is $\dps{\sod X}{\s_d}$-nonmeager ($\dps{\sod X}{\s_d}$-Baire);
\item[\rm (2)] $\dd h$-nonmeager ($\dd h$-Baire) if $X$ is $\dps{\sod X}{\s_h}$-nonmeager ($\dps{ \sod X}{\s_h}$-Baire);
\item[\rm (3)] $\dd s$-nonmeager ($\dd s$-Baire) if $X$ is $\dps{\sod X}{\s_s}$-nonmeager ($\dps{ \sod X}{\s_s}$-Baire);
\item[\rm (4)] $\rd \ga$-nonmeager ($\rd \ga$-Baire) if $X$ is $\dps{\pp_\ga(\sod X)}{\s_d} $-nonmeager ($\dps{\pp_\ga(\sod X)}{\s_d}$-Baire);
\item[\rm (5)] $\lrd c\ga$-nonmeager ($\lrd c\ga$-Baire) if $X$ is $\dps{\pp^c_\ga(\sod X)} {\s_d}$-nonmeager ($\dps{\pp^c_\ga(\sod X)}{\s_d}$-Baire);
\end{itemize}
\end{definition}

\begin{note} The spaces that are here called "$\D$-nonmeager"\ were called "$\D$-Baire" in papers \cite{rezn2008,tk2014,moors2017}.
\end{note}

\begin{proposition}[Proposition \dref{pdp1} \cite{rezn2022gbd}] \label{pdp1}
Let $X$ be a space, $\s$ be \nfs/ and $n\in\om$.
Consider the condition
\begin{itemize}
\item[\rm(*)]
\begin{enumerate}
\item
A space $X$ is $\tD(\s)$-nonmeager if and only if some nonempty open $U\subset X$ is $\tD(\s)$-Baire.
\item
If the space $X$ is homogeneous, then $X$ is $\tD(\s)$-nonmeager if and only if $X$ is $\tD(\s)$-Baire.
\end{enumerate}
\end{itemize}
\begin{itemize}
\item[\rm (1)] If $\tD\in\sset{\D,\rd n}$ then {\rm (*)} holds.
\item[\rm (2)] If the space $X$ is quasiregular and $\tD\in\sset{\D,\rd n, \ld v, \ld c, \lrd c+, \lrd cn }$, then {\rm (*)} holds.
\end{itemize}
\end{proposition}

\begin{theorem}[Theorem \dref{tgd1} \cite{rezn2022gbd}] \label{tgd1}
Let $\cP_k$ ($\cP_c$) be the smallest class of spaces that
\begin{itemize}
\item contains $p$-spaces and strongly $\Sigma$-spaces;
\item is closed under arbitrary (countable) products;
\item is closed under taking open subspaces.
\end{itemize}

For $\wt\D\in \sset{\D_s,\D_h,\D}$, if a regular Baire (nonmeager) space $X$ belongs to the class of spaces described in $(\wt\D)$, then $ X$ is $\wt\D$-Baire ($\wt\D$-nonmeager).
\begin{itemize}
\item[$(\D_s)$] $\sigma$-spaces.
\item[$(\D_h)$] is the smallest class of spaces that
\begin{itemize}
\item contains $\Sigma$-spaces and $w\D$-spaces;
\item is closed under products by spaces from the class $\cP_c$;
\item is closed under taking open subspaces.
\end{itemize}
\item[$(\D)$] is the smallest class of spaces that
\begin{itemize}
\item contains $\Sigma$-spaces, $w\D$-spaces, and feebly compact spaces;
\item is closed under products by spaces from the class $\cP_k$;
\item is closed under taking open subspaces.
\end{itemize}
\end{itemize}
\end{theorem}

\subsection{\cdp/-Baire spaces}\label{sec-cdp}

Let $X$ be a space, $\cG$ be a family of open subsets of $X$.
Denote
\begin{align*}
B(X,\cG) \eqdef \{ x\in X:\ & \text{ family }\\
&\set{\st(x,\ga):\ga\in\cG\text{ and }x\in\bigcup\ga}
\\
&\ \text{ is the base at point }x \}
\\
\dev X\eqdef \min \{\, |\cG| :\
&\cG \text{ open cover family }X
\\
& \text{ and }X=B(X,\cG)
,\}.
\end{align*}

Let $Y\subset X$, $\cN$ be the family of all nowhere dense subsets of $X$. If $Y\not \subset \bigcup \cN$ then put $\nov {Y,X} \eqdef \infty$, otherwise
\begin{align*}
\nov {Y,X} &\eqdef \min \set{|\cL|: \cL\subset \cN\text{ and } Y\subset \bigcup\cL},
\\
\nov {X} &\eqdef \nov {X,X}.
\end{align*}

\begin{definition}
Let $X$ be a space.
\begin{itemize}
\item
We call a space a \term{\cdp/-space} if $\dev X<\nov X$, that is, if there exists a family $\cP$ of open coverings of the space $X$such that $|\cP|<\nov X $ and $B(X,\cP)=X$.
\item
We call a space a \term{\cdpz/-space} if there exists a family $\cP$ of open partitions of $X$ such that $|\cP|<\nov X$ and $B(X,\cP)=X$ .
\item We call a space $X$ \term{\cdp/-nonmeager} if there exists a family $\cG$ of open families in $X$ such that $|\cG|<\nov{Y,X}$ for $Y =B(X,\cP)$.
\item We call a space $X$ \term{\cdp/-Baire} if
every nonempty open subset of $X$ is a \cdp/-nonmeager space.
\end{itemize}
\end{definition}

In the article \cite{cdp2010} one studies spaces for which $\dev X<\nov X$, that is, \cdp/-spaces. Metrizable nonmeager spaces are \cdp/-spaces. Examples of non-metrizable \cdp/-spaces can be obtained using Martin's axiom ($MA$).

\begin{proposition}[Proposition \dref{pcdp2} \cite{rezn2022gbd}] \label{pcdp2}
($MA$)
Let $\tau<2^\om$ be an infinite cardinal, $X$ be an absolute $G_\tau$ space with countable Suslin number and $\dev X\leq \tau$, for example $X=\R^\tau$. Then $X$ is a \cdp/-space.
\end{proposition}

\begin{cor} \label{ccdp1} \rm{(MA)}
Let $X$ be a ccc compact space and $w(X)<2^\om$.
Then $X$ is a \cdp/-space.
\end{cor}

Clearly, a \cdp/-space is a \cdpz/-space.

\begin{theorem}[Theorem \dref{tcdp1} \cite{rezn2022gbd}] \label{tcdp1}
Let $X$ be a semiregular space. The following conditions are equivalent
\begin{itemize}
\item[\rm (1)] $X$ is a \cdp/-nonmeager space;
\item[\rm (2)] some nonempty open subset $U\subset X$ contains a dense \cdp/-space $Y\subset U\subset \cl Y$;
\item[\rm (3)] some nonempty open subset $U\subset X$ contains a dense \cdpz/-space $Y\subset U\subset \cl Y$;
\end{itemize}
\end{theorem}

Theorem \ref{tcdp1} implies

\begin{proposition}[Proposition \dref{pcdp5} \cite{rezn2022gbd}] \label{pcdp5}
Let $X$ be a semiregular space.
If $X$ has a metrizable nonmeager subspace, then $X$ is a \cdp/-nonmeager space.
\end{proposition}

\subsection{\texorpdfstring{$\D_s$}{Ds}-Baire spaces}\label{sec-dsp}

Known examples of $\D_s$-Baire spaces are obtained using the following proposition

\begin{proposition}[Proposition \dref{pdsp1} \cite{rezn2022gbd}] \label{pdsp1}
Let $X$ be a space. If $X$ is \cdp/-nonmeager (\cdp/-Baire), then $X$ is $\D_s$-nonmeager ($\D_s$-Baire).
\end{proposition}

\cite{rezn2022gbd} lists other subclasses of $\D_s$-nonmeager spaces, but it is not known whether these classes contain non-\cdp/-nonmeager spaces.

Propositions \ref{pcdp5} and \ref{pdsp1} imply

\begin{proposition}\label{pdsp2} 
Let $X$ be a space. If $X$ contains a nonmeager metrizable subspace, then $X$ is $\D_s$-nonmeager.
\end{proposition}

Under the assumption MA+$\lnot$CH, there exists a $\D_s$-nonmeager space without first-countable points:
Propositions \ref{pdsp1} and \ref{pcdp2} imply that the space $\R^{\om_1}$ is $\D_s$-nonmeager.

\begin{proposition}\label{pdsp3} 
Let $(X,\tp)$ be a space, $\tps=\tp\setminus \sset\es$. The following conditions are equivalent:
\begin{itemize}
\item[\rm (1)] $X$ is a $\D_s$-nonmeager space.
\item[\rm (2)] Let $W:\tps\to\tps$ be a mapping such that
\[
W(V)\subset W(U)\subset U
\]
for $V,U\in\tps$, $V\subset U$. Then there exists $x\in X$ such that $x\in W(U)$ for every $U\in\St (x,\tps)$.
\end{itemize}
\end{proposition}
\begin{proof}
(1)$\implies$(2). Let us assume the opposite. Then there exists a semi-open neighborhood $P$ of the diagonal such that $x\notin W(P(x))$ for all $x\in X$. Since $X$ is a $\D_s$-nonmeager space, there exists $U\in\tps$ for which $M=\set{x\in U: U\subset P(x)}$ is dense in $U $. Let $x\in W(U) \cap M$. Then $x\in U\subset W(P(x))$. Contradiction.

(2)$\implies$(1). Let us assume the opposite. Then there exists a semi-open diagonal neighborhood $P$ such that $M(U)=\set{x\in U: U\subset P(x)}$ is not dense in $U$ for all $U\in\tps$ . Let $W(U)=U\setminus \cl{M(U)}$. Conditions (2) are satisfied for $W$. Then there exists $x\in X$ for which $x\in W(U)$ for every $U\in\St (x,\tps)$. Then $x\in W(U)$ for $U=P(x)$. Hence $x\notin M(U)$ and $U\not\subset P(x)$. A contradiction, since $U=P(x)$.
\end{proof}

\begin{theorem} \label{tdsp2} 
Let $X$ be a $\D_s$-nonmeager space and let $f: X\to X$ be a feebly homeomorphism. Then for some $x\in X$ the mapping $f$ is continuous in $x$ and $f^{-1}$ is continuous in $f(x)$.
\end{theorem}
\begin{proof} Let $\tps$ be all open nonempty subsets of $X$. We set
\[
W(U)=\Int f^{-1}(\Int f(U))
\]
for $U\in\tps$. Since $f$ is a feebly homeomorphism, $W(U)\neq\es$ and $W$ satisfy the conditions of Proposition \ref{pdsp3}.
Then there exists $x\in X$ such that $x\in W(U)$ for every $U\in\St (x,\tps)$. The point $x$ is the desired one.
\end{proof}

\section{Properties of groups defined by invariant semi-neighborhoods of the diagonal}
\label{sec-rtgp}
In this section, $(G,\tp)$ is a right semitopological group, $\tps=\tp\setminus\sset\es$, $\nne$ is a family of open neighborhoods of the identity, $\nnee$ is a family of neighborhoods of the identity that are not necessarily open.

\subsection{Diagonal Semi-Neighborhoods in Groups}\label{sec-gsonbd}

A set $P\subset G\times G$ is called \term{right-invariant} if
\[
\set{(xg,yg): (x,y)\in P}=P
\]
for all $g\in G$. For $M\subset G$ we set
\[
\nn(M) \eqdef \set{(x,gx): g\in M\text{ and }x\in G}.
\]
If $P=\nn(M)$, then $P(g)=Mg$.
The set $\nn(M)$ is right-invariant. The mapping $\nn$ establishes a one-to-one correspondence between subsets of $G$ and right-invariant subsets of $G\times G$. If $P\subset G\times G$ is a right-invariant subset, then $P=\nn(P(e))$.
Denote
\begin{align*}
\clg M &\eqdef \bigcap\set{\cl{MV}: V\in\nne },
\\
\clh M &\eqdef \bigcap\set{MV: V\in\nne }.
\end{align*}

\begin{assertion} \label{agsonbd1}
Let $L,M\subset G$. Then
\begin{itemize}
\item[\rm (1)] $(x,y)\in \nn(M)$ if and only if $yx^{-1}\in M$;
\item[\rm (2)] $(\nn(M))^{-1}=\nn(M^{-1})$;
\item[\rm (3)] $\cl{\nn(M)}=\nn(\clg M)$;
\item[\rm (4)] $\nn(M)\circ \nn(L)=\nn(LM)$.
\end{itemize}
\end{assertion}
\begin{proof}
(1) $(x,y)\in \nn(M)$ $\iff$ $y=gx$ for some $g\in M$ $\iff$ $yx^{-1}\in M$.

(2) $(x,y)\in (\nn(M))^{-1}$ $\iff$ $(y,x)\in \nn(M)$ $\iff$ $xy^{ -1}\in M$ $\iff$ $(x,y)\in \nn(M^{-1}$.

(3) $(x,y)\in \cl{\nn(M)}$ $\iff$ $(Vx\times Uy) \cap \nn(M)\neq\es$ for any $V,U \in\nne$
$\iff$ $MUx \cap Vy\neq\es$ for any $V,U\in\nne$ $\iff$ $y\in \cl{MUx}$ for any $U\in\nne$ $\ iff$ $y\in \clg{M}x$.

(4) $(x,y)\in \nn(M)\circ \nn(L)$ $\iff$ $(x,z)\in \nn(M)$ and $(z,y)\in  \nn(L)$ for some $z\in G$ $\iff$ $z\in Mx$ and $y\in Lz$ for some $z\in G$ $\iff$ $y\in LMx$ .
\end{proof}

Denote
\[
\gsod G\eqdef \set{P\in\sod G: G\text{ is a right invariant subset}}.
\]
The mapping $\nn$ establishes a bijection between $\nnee$ and $\gsod G$.

On $\Exp{\nnee }$ we introduce an order relation and an equivalence relation related to it. For $\gcP,\gcQ\subset \nnee $ we put
\begin{itemize}
\item
$\gcP\lleq \gcQ$ if and only if for any $V\in \gcQ$ there exists $U\in\gcP$such that $U\subset V$;
\item
$\gcP\sim \gcQ$ if and only if $\gcP\lleq \gcQ$ and $\gcQ\lleq \gcP$.
\end{itemize}
Note that if $\gcP\supset \gcQ$, then $\gcP\lleq \gcQ$.
For $\gcP\subset \nnee $ we denote
\[
\gpp^e(\gcP)\eqdef\set{W\in \nnee : U\subset W,\text{ for some }U\in\gcP}.
\]

\begin{proposition} \label{gpsonbd2}
For $\gcP,\gcQ\subset\nnee$
\begin{itemize}
\item
$\gcP\lleq \gcQ$ if and only if $\gpp^e(\gcP)\supset \gpp^e(\gcQ)$;
\item
$\gcP\sim \gcQ$ if and only if $\gpp^e(\gcP) = \gpp^e(\gcQ)$.
\end{itemize}
\end{proposition}

For $\gcP\subset \nnee $ we denote
\begin{align*}
\gpp(\gcP)&\eqdef\gcP,
\\
\gpp_+(\gcP)&\eqdef\set{UV: V\in\nnee \text{ and }U\in \gcP},
\\
\gpp^v(\gcP)&\eqdef\set{\cl U: U\in\gcP},
\\
\gpp^c(\gcP)&\eqdef\set{\clg U: U\in\gcP},
\\
\gpp_+^c(\gcP)&\eqdef \gpp^v(\gpp_+(\gcP)).
\end{align*}

It follows from the definitions

\begin{proposition} \label{gpsonbd3}
Let $X$ be a space.
For $\gcP\subset\nnee$
\begin{itemize}
\item
$\gcP\lleq \gpp_+(\gcP)\lleq \gpp_+^c(\gcP)$
\item
$\gcP\lleq \gpp^v(\gcP)\lleq \gpp^c(\gcP)\lleq \gpp_+^c(\gcP)$
\end{itemize}
\end{proposition}

Let $\la$ be an ordinal, $M_\al\subset G$ for $\al<\la$. By induction on $\be<\ga$ we define the sets $\gppl_\be(\cp_\be),\gppl^c_\be(\cp_\be)\subset G$, where $\cp_\be=(M_\al)_{\al<\be}$.

$\be=0$. We set
\[
\gppl_0(\cp_0)\eqdef\gppl_0^c(\cp_0)\eqdef\sset\es.
\]

$\be=1$. We set
\[
\gppl_1(\cp_1)\eqdef M_0,
\qquad
\gppl_1^c(\cp_1)\eqdef \clg{M_0}.
\]

$1<\be\leq \ga$. We set
\begin{align*}
\gppl_\be(\cp_\be)&\eqdef\begin{cases}
Q_\be,&\text{if $\be$ is cardinal limit}
\\
M_{\be'} Q_\be,&\text{if $\be=\be'+1$}
\end{cases}
\\
\gppl^c_\be(\cp_\be)&\eqdef\begin{cases}
\clg{Q^c_\be},&\text{if $\be$ is cardinal limit}
\\
\clg{M_{\be'} Q^c_\be},&\text{if $\be=\be'+1$}
\end{cases}
\end{align*}
where
\[
Q_\be=\bigcup_{\al<\be} \gppl_\al(\cp_\al)
\qquad
\text{and}
\qquad
Q^c_\be=\bigcup_{\al<\be} \gppl^c_\al(\cp_\al).
\]

For $\gcP\subset \nnee $ we denote
\begin{align*}
\gpp_\ga(\gcP)&\eqdef\set{\gppl_\ga(\cp): \cp\in \gcP^\ga},
\\
\gpp^c_\ga(\gcP)&\eqdef\set{\gppl^c_\ga(\cp): \cp\in \gcP^\ga}.
\end{align*}

From the clause \ref{gpsonbd3} and the definitions it follows


\begin{proposition} \label{pgsonbd4}
Let $\gcP \subset \nnee $ and $1<\al<\be$ be ordinals.
Then
\begin{align*}
\gpp_+(\gcP) \lleq \gpp_2(\gcP) \lleq &\gpp_\al(\gcP) \lleq \gpp_\be(\gcP) \lleq \gpp_\be^c(\gcP)
\\
\gpp^c(\gcP) = \gpp_1^c(\gcP) \lleq &\gpp_\al^c(\gcP) \lleq \gpp_\be^c(\gcP)
\\
\gpp^c(\gcP) \lleq \gpp^c_+(\gcP) \lleq &\gpp_2^c(\gcP)
\end{align*}
\end{proposition}

For $\cB\subset \nnee$, $\cB\sim \nnee$ if and only if $\cB$ is a base in $e$.

\begin{assertion} \label{agsonbd2}
Let $G$ be an \rr-paratopological group.
\begin{itemize}
\item[\rm (1)] $\gpp_\ga(\nnee)\sim \nnee$ for $0<\ga\leq \om$.
\item[\rm (2)] If $G$ is a regular space, then $\gpp_\ga^c(\nnee)\sim \nnee$ for $0<\ga<\om$.
\end{itemize}
\end{assertion}
\begin{proof} (1) Let $U\in\nnee$. There is $\cp=\sq{U_n}{n<\la}\subset \nnee$, so $U_0^2\subset U$ and $U_{n+1}^2\subset U_n$ for $n< \la$. Then $V=\gppl_\ga(\cp)\in \gpp_\ga(\nnee)$ and $V\subset U$.

(2) Let us prove by induction on $\ga$. For $\ga=1$ the assertion is obvious. Let $\ga=n+1$ and $U\in\nnee$. Then $\clg{V^2}\subset \cl{V^3} \subset U$ for some $V\in\nnee$. By the inductive hypothesis $S\subset V$ for some $S\in \gpp^c_n(\nnee)$. Then $Q=\clg{V S}\in \gpp^c_\ga(\nnee)$ and $Q\subset U$.
\end{proof}

From the definitions and Statement \ref{agsonbd1} it follows
\begin{proposition} \label{pgsonbd5}
Let $\gcP \subset \nnee$, $\cP=\nn(\gcP)$ and $\ga$ be an ordinal.
Then
\begin{align*}
\pp^v(\cP)&=\nn(\gpp^v(\gcP)),
&
\pp^c(\cP)&=\nn(\gpp^c(\gcP)),
\\
\pp_+(\cP)&\lleq\nn(\gpp_+(\gcP)),
&
\pp_+^c(\cP)&\lleq\nn(\gpp_+^c(\gcP)),
\\
\pp_\ga(\cP)&=\nn(\gpp_\ga(\gcP)),
&
\pp_\ga^c(\cP)&=\nn(\gpp_\ga^c(\gcP)).
\end{align*}
\end{proposition}

\subsection{\texorpdfstring{$\gdps \gcP\s$}{gDNS}-Baire groups}\label{sec-gdps}

\begin{definition}
Let $\s$ be the normal square functor, $\gcP\subset \nnee$. Let us call the group $G$
\term{$\gdpss \gcP$-Baire} ($\gdpss \gcP$-Baire) group if $G$ is a $\dpss{\nn(\gcP)}$-Baire space.
\end{definition}

Any element of $P\in \nn(\gcP)$ is a right-invariant subset and $\rho_g\times \rho_g(\sorestr{P}W) = \sorestr{P}{\rho_g(W)}$ for $g\in  G$ and $W\in\tps$. So if $\sorestr{P}W\in\s(W)$ then $\sorestr{P}{Wg}\in\s(Wg)$.
We get that $G$ is a $\dpss{\nn(\gcP)}$-Baire space if and only if $G$ is a $\dpss{\nn(\gcP)}$-nonmeager space.

It follows from the definitions

\begin{proposition} \label{pgdps1}
Let $\s$, $\gr$ be \nfs/, $\gcP,\gcQ\subset \nnee$.
Suppose $\s\ggeq \gr$ and $\gcP\lleq\gcQ$.
If $G$ is a $\gdps \gcP\s$-Baire group, then $G$ $\gdps \gcQ\gr$ is a Baire group.
\end{proposition}

Similar to the Proposition \ref{pdps7}, it checks

\begin{proposition} \label{pgdps7}
Let $\gcP \subset \nnee$. If in the diagram from Proposition \ref{pdps7} replace $\pp$ with $\gpp$, then
arrow
\[
\ddpsV Fk \to \ddpsV Gl
\]
means that $F,G:\Exp{\nnee}\to \Exp{\nnee}$ are mappings, $k,l\in\sset{d,v,h,s,a}$ and the following condition is satisfied:
\begin{itemize}
\item[]
if $G$ is a $\gdps{F(\gcP)}{\s_k}$-Baire group, then $G$ is a $\gdps{G(\gcP)}{\s_l}$-Baire group.
\end{itemize}
\end{proposition}

\begin{proposition} \label{pgdps8-1}
Let $G$ be an \rr-semitopological group, $M\subset G$. Then
\[
\clh M = \cl{M^{-1}}^{\,-1}.
\]
\end{proposition}
\begin{proof} Let $x\in G$. Then $x\in\clh M$ if and only if for any $V\in\nne$ $x\in MV$
$\iff$ $x^{-1}\in V^{-1}M^{-1}$ $\iff$ $Vx^{-1}\cap M^{-1}\neq\es$ , that is, $x^{-1} \in\cl{M^{-1}}$ and $x \in\cl{M^{-1}}^{-1}$.
\end{proof}

\begin{proposition} \label{pgdps8}
Let $U\in \nnee$ and $W\in\tps$. For $k\in\sset{d,v,h,s,a}$, $\sorestr{\nn(U)}{W}\in\s_k(W)$ if and only if the condition \cnLg k :
\begin{itemize}
\item[\cnLg d] $WW^{-1}\subset \clg{U}$;
\item[\cnLg v] $WW^{-1}\subset \cl{U}$;
\item[\cnLg h] $WW^{-1}\subset \cl{U^{-1}}$ or equivalently $WW^{-1}\subset \clh{U}$;
\item[\cnLg s] $WM^{-1}\subset U$ for some $M\subset W\subset \cl M$;
\item[\cnLg a] $WW^{-1}\subset U$.
\end{itemize}
\end{proposition}
\begin{proof}
\begin{itemize}
\item[\cnLg a]
For $x,y\in W$,
$(x,y)\in\nn(U)$ if and only if $yx^{-1}\in U$.
\item[\cnLg s]
For $x\in M$,
$\sset x \times W \subset \nn(U)$ if and only if $Wx^{-1}\in U$.
\item[\cnLg h]
$W\subset \cl{U^{-1}x}$ for all $x\in W$ if and only if $WW^{-1}\subset \cl{U^{-1}}$. The Proposition \ref{pgdps8-1} implies that $WW^{-1}\subset \cl{U^{-1}}$ is equivalent to $WW^{-1}\subset \clh{U} $.
\item[\cnLg v]
$W\subset \cl{Ux}$ for all $x\in W$ if and only if $WW^{-1}\subset \cl{U}$.
\item[\cnLg v]
Follows from $\cl{\nn(U)}=\nn(\clg U)$ and \cnLg a.
\end{itemize}
\end{proof}

\subsection{\texorpdfstring{$\gds \s$}{gDS}-Baire groups}\label{sec-gdp}

\begin{definition}
Let $\s$ be a normal square functor, $G$ be an \rr-topological group, $\ga$ be an ordinal,
$\gcP = \nnee$. Let $\gpp^k_l: \Exp{\nnee }\to \Exp{\nnee }$ be one of the mappings considered in the \ref{sec-gsonbd} section:
\[
\gpp^k_l\in\sset{\gpp,\gpp^e,\gpp^v,\gpp^c,\gpp_+,\gpp_+^c,\gpp_\ga,\gpp_\ga^c}.
\]
We say $G$ is a \term{$\glrds kl{\s}$-Baire} group if $G$ is a $\gdps{\gpp^k_l(\nnee )}\s$-Baire group. For
\[
\tD\in\sset{\gD,\gld e,\gld v,\gld c,\grd +,\glrd c+,\grd\ga,\glrd c\ga}
\]
we have defined $\tD(\s)$-Baire groups. For
\[
k\in\sset{d,h,v,s,a}
\]
we say $G$ is a \term{${\tD}_k$-Baire} group if $G$ is a $\tD(\s_k)$-Baire group. For
\[
\til{\tD}\in\sset{\gdd k, \gldd {e}k,\gldd {v}k,\gldd {c}k,\grdd {+}k,\glrdd c{_+} k,\grdd {_\ga}k,\glrdd c{_\ga}k}
\]
we have defined $\til{\tD}$-Baire groups. Also, if the subscript is not written, $d$ is implied: $\glrd kl$-tober is $\glrdd kld$-tober. A separate direct definition for the most important classes of groups: the group $G$ is called
\begin{itemize}
\item
$\gD$-Baire if $G$ is $\gdps{\nnee }{\s_d}$-Baire;
\item
$\gdd h$-Baire if $G$ is $\gdps{\nnee }{\s_h}$-Baire;
\item
$\gdd s$-Baire if $G$ is $\gdps{\nnee }{\s_s}$-Baire;
\item
$\gdd v$-Baire if $G$ is $\gdps{\nnee }{\s_v}$-Baire;
\item
$\gdd a$-Baire if $G$ is $\gdps{\nnee }{\s_a}$-Baire;
\item
$\grd \ga$-Baire if $G$ is $\gdps{\gpp_\ga(\nnee )}{\s_d}$-Baire;
\item
$\glrd c\ga$-Baire if $G$ is $\gdps{\gpp^c_\ga(\nnee )}{\s_d}$-Baire.
\item
$\grdda\ga$-Baire if $G$ is $\gdps{\gpp_\ga(\nnee )}{\s_a}$-Baire;
\end{itemize}
\end{definition}

Since $\gpp_+(\nnee )=\gpp_2(\nnee )$, the classes of $\grd +(\s)$-Baire and $\grd 2(\s)$-Baire groups coincide.

Since $\nn(\nnee)=\gsod G\subset \sod G$, Proposition \ref{pgsonbd5} implies

\begin{proposition} \label{pgdp1-1}
Let $\s$ be a normal square functor, $G$ be an \rr-topological group, $\ga$ be an ordinal,
\[
\gpp^l_m\in\sset{\gpp,\gpp^e,\gpp^v,\gpp^c,\gpp_+,\gpp_+^c,\gpp_\ga,\gpp_\ga^c}.
\]
If $G$ is an $\lrds lm{\s}$-Baire space, then $G$ is an $\glrds lm{\s}$-Baire group.

Let
\[
\glrdd l{_m}k\in\sset{\gdd k, \gldd {e}k,\gldd {v}k,\gldd {c}k,\grdd {+}k,\glrdd c{_+ }k,\grdd {_\ga}k,\glrdd c{_\ga}k}
\]
where $k\in\sset{d,h,v,s,a}$. If $G$ is an $\lrdd l{_m}k$-Baire space, then $G$ is a $\glrdd l{_m}k$-Baire group.
\end{proposition}

From Proposition \ref{pgdps8} it follows

\begin{proposition} \label{pgdp1}
Let
\[
\tD\in\sset{\gD,\gdd h,\gdd s, \gdd v, \gdd a}.
\]
A group $G$ is $\tD$-Baire if for any $U\in\nne$ there exists $W\in\nne$ such that the condition $(\tD)$ is satisfied:
\begin{itemize}
\item[($\gD$)] $WW^{-1}\subset \clg{U}$;
\item[($\gdd h$)] $WW^{-1}\subset \cl{U^{-1}}$ or equivalently $WW^{-1}\subset \clh{U}$ ;
\item[($\gdd s$)] $WM^{-1}\subset U$ for some $M\subset W\subset \cl M$;
\item[($\gdd v$)] $WW^{-1}\subset \cl{U}$;
\item[($\gdd a$)] $WW^{-1}\subset U$.
\end{itemize}
\end{proposition}

\begin{proposition} \label{pgdp2}
Let $\ga$ be an ordinal and
\[
\tD\in\sset{\grd \ga,\glrd c\ga,\grdda\ga}.
\]
A group $G$ is $\tD$-Baire if for any $\cp=\sq{W_\al}{\al<\ga}\in \nne^\ga$ there exists $W\in\nne$, so the condition $(\tD)$ is satisfied:
\begin{itemize}
\item[($\grd \ga$)] $WW^{-1}\subset \clg{\gppl_\ga(\cp)}$;
\item[($\glrd c\ga$)] $WW^{-1}\subset \gppl_\ga^c(\cp)$;
\item[($\grdda\ga$)] $WW^{-1}\subset \gppl_\ga(\cp)$.
\end{itemize}
\end{proposition}

From Proposition \ref{pdps7} it follows

\begin{proposition} \label{pgdp3}
In the diagrams below, the arrow
\[
A \to B
\]
means that
if $G$ is an $A$-Baire group, then $G$ is a $B$-Baire group.
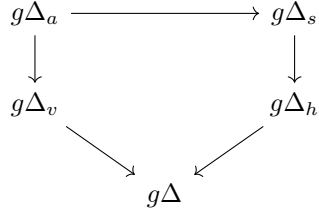
\[
\begin{tikzcd}
\gdd a \ar[d]\ar[rr]&&\gdd s \ar[d]
\\
\gdd v \ar[dr] &&\gdd h \ar[dl]
\\
&\gD&
\end{tikzcd}
\]
\end{proposition}

\section{Continuity in \rr-semitopological groups}\label{sec-rtgc}

\begin{theorem}\label{tgtp1}
\hspace{2em}
\begin{itemize}
\item[{\rm (1)}]
A \rr-semitopological group $G$ is an \rr-topological group if and only if $G$ is $\gdd a$-Baire.
\item[{\rm (2)}]
A \rr-quasitopological group $G$ is $\gdd v$-Baire if and only if $G$ is $\gdd h$-Baire.
\item[{\rm (3)}]
A semiregular \rr-semitopological group $G$ is an \rr-topological group if and only if $G$ is $\gdd v$-Baire.
\item[{\rm (4)}]
A semiregular \rr-quasitopological group $G$ is an \rr-topological group if and only if $G$ is $\gdd h$-Baire.
\item[{\rm (5)}]
A semiregular \rr-paratopological group $G$ is an \rr-topological group if and only if $G$ is $\gD$-Baire.
\item[{\rm (6)}]
Let $\ga\leq\om$.
A \rr-paratopological group $G$ is an \rr-topological group if and only if $G$ is $\grdda\ga$-Baire.
\item[{\rm (7)}]
Let $\ga<\om$.
A regular \rr-paratopological group $G$ is an \rr-topological group if and only if $G$ is $\glrd c\ga$-Baire.
\end{itemize}
\end{theorem}
\begin{proof} (1) Follows from Proposition \ref{pgdp1} ($\gdd a$).

\noindent
(2) Follows from Proposition \ref{pgdp1} ($\gdd v$) and ($\gdd h$).

\noindent
(3) Let $U\in\nne$ be a regular open set. Proposition \ref{pgdp1} ($\gdd v$) implies that $WW^{-1}\subset \cl U$ for some $W\in\nne$. Since $WW^{-1}$ is open, $WW^{-1}\subset \Int\cl U=U$.

\noindent (4)
Follows from (2) and (3).

\noindent (5)
Let $U\in\nne$ be a regular open set. For some $V\in \nne$, $V^2\subset U$. Proposition \ref{pgdp1} ($\gD$) implies that $WW^{-1}\subset \clg V\subset \cl{V^2}\subset \cl U$ for some $W\in\nne$.
Then $WW^{-1}\subset \Int \cl U= U$.

\noindent (6)
Follows from Proposition \ref{pgdp2} ($\grdda\ga$) and Statement \ref{agsonbd2} (1).

\noindent (7)
Follows from Proposition \ref{pgdp2} ($\glrd c\ga$) and Statement \ref{agsonbd2} (2).
\end{proof}

\begin{theorem}\label{tgtp2}
Let $G$ be an \rr-semitopological $\gdd s$-Baire group.
Then $\tcs_f(G)$ is a topological group and $\la_g$ is a homeomorphism for all $g\in \tcs_f(G)$.
\end{theorem}
\begin{proof}
The set $H=\tcs_f(G)$ is a subgroup of $G$.
Let $g\in H$. Let us show that $\la_g$ is a homeomorphism. For $U\in\tps$ we denote
\[
\wt U = \Int g^{-1}\Int gU = \Int \la_g^{-1}(\Int \la_g(U)).
\]
Note that $\wt U\in\tps$, $\wt U\subset U\subset \cl{\wt U}$, $\wt{Ux}= \wt{U}x$ for $x\in G$ and $\wt V\subset \wt W$ for $V\subset U$, $V\in\tps$.
Assume that $\la_g$ is not a homeomorphism. Then $x\notin \wt{Ux}=\wt{U}x$ for some $x\in G$ and $U\in\tps$. Then $e\notin \wt U$.
\begin{lemma}\label{lgtp0}
Let $U\in\nne$. Then $WM^{-1}\subset U$ for some $W\in\tps$ and $M\subset W\subset \cl M$ so $e\in M$ and $W\subset U$.
\end{lemma}
\begin{proof} Since $G$ is a $\gdd s$-Baire group, it follows from Proposition \ref{pgdps8} that $W_*M_*^{-1}\subset U$ for some $W_* \in \nne$ and $M_*\subset W\subset \cl M_*$. Let $W=W_*\cap U$ and $M=M_*\cup \sset e$.
\end{proof}
Lemma \ref{lgtp0} implies that $WM^{-1}\subset U$ for some $W\in \nne$, $W\subset U$ and $M\subset W\subset \cl M$. Since $\wt W$ is open and dense in $W$ and $M$ is dense in $W$, there exists $y\in \wt W\cap M$. Since $Wy^{-1}\subset U$, then
\[
e\in \wt W y^{-1}\subset \wt {W y^{-1}} \subset \wt U.
\]
A contradiction with the fact that $e\notin \wt U$.

The group $H$ is a semitopological group.
\begin{lemma}\label{lgtp1}
Let $M\subset G$, $W\in\nne$. Then
$H\cap \cl M \subset MW^{-1}$.
\end{lemma}
\begin{proof}
Let $x\in H\cap \cl M$. Then $xW\cap M\neq\es$. Hence $x\in MW^{-1}$.
\end{proof}
Let us show that $H$ is a quasitopological group, that is, for $U\in\nne$ there exists $W\in\nne$ such that $H\cap W\subset U^{-1}$. Lemma \ref{lgtp0} implies that $WM^{-1}\subset U$ for some $W\in \nne$ and $M\subset W\subset \cl M$. Lemma \ref{lgtp1} implies that
\[
H\cap W \subset H\cap \cl M \subset MW^{-1} \subset U^{-1}.
\]

Let us show that $H$ is a paratopological group, that is, for $U\in\nne$ there exists an $S\in\nne$ such that $(H\cap S)^2\subset U$. Lemma \ref{lgtp0} implies that $WM^{-1}\subset U$ for some $W\in \nne$ and $M\subset W\subset \cl M$ such that $e\in M$ and $W\subset U$. Also $VL^{-1}\subset W$ for some $V\in \nne$ and $L\subset V\subset \cl V$ so $e\in L$ and $V\subset W$. Then $MLV^{-1} \subset U^{-1}$. Lemma \ref{lgtp1} implies that $H\cap \cl{ML}\subset MLV^{-1} \subset U^{-1}$. Let $W_*=W\cap H$ and $V_*=V\cap H$. Since $\cl M L\subset \cl{ML}$ and $W_*\subset \cl M$, then $W_*L\subset \cl{ML}$. Since $V_*\subset \cl L$ and $W_*\subset H$, then $W_*V_*\subset W_*\cl L\subset \cl{ML}$. We get $W_*V_*\subset U^1$. Since $V_*\subset W_*$, then $(V\cap H)^2\subset U^{-1}$. Since $H$ is a quasitopological group, $S\cap H\subset V^{-1}$ for some $S\in\nne$. Then $(H\cap S)^2\subset U$.
\end{proof}

Theorem \ref{tgtp2} implies
\begin{theorem}\label{tgtp3}
Let $G$ be a \rr-semitopological \goo{fl}-topological $\gdd s$-Baire group.
Then $G$ is a topological group.
\end{theorem}

\begin{lemma}\label{lgtp3}
Let $G$ be an \rr-topological group, $M\subset G$, $W\in\nne$. Then $\cl M \subset W^{-1}M$.
\end{lemma}
\begin{proof}
Let $x\in \cl M$. Then $Wx\cap M\neq\es$. Therefore, $x\in W^{-1}M$.
\end{proof}

\begin{theorem}\label{tgtp4}
Let $G$ be a \rr-semitopological \goo{dfl^*}-topological $\gdd s$-Baire group.
Then
\begin{itemize}
\item[{\rm (1)}]
$G$ is a $\gdd v$-Baire group;
\item[{\rm (2)}]
if $G$ is semiregular, then $G$ is a topological group.
\end{itemize}
\end{theorem}
\begin{proof}
Condition (2) follows from condition (1), Theorem \ref{tgtp1} (3) and Theorem \ref{trtg2}. Let us prove (1). 
Theorem \ref{tgtp2} implies that $H=\tcs(G)$ is dense in $G$ and $H$ is a topological group.
\begin{lemma}\label{lgtp2}
Let $U\in\nne$. There exists $V\in\nne$ such that $\cl{V}^2\subset \cl U$.
\end{lemma}
\begin{proof}
There is $V\in\nne$ such that $V_*=V_*^{-1}$ and $V_*^2\subset U$, where $V_*=V\cap H$. Then $V_*\cl{V_*}\subset \cl U$ and $\cl{V_*}\,\cl{V_*}\subset \cl U$. Since $\cl{V_*}=\cl{V}$, then $\cl V^2\subset \cl U$.
\end{proof}
By Proposition \ref{pgdp1} ($\gdd s$), it suffices to show that for $U\in\nne$ there exists a $W\in\nne$ such that $WW^{-1}\subset \cl U$. Lemma \ref{lgtp2} implies that there exist $V,S\in \nne$ such that $\cl V^2\subset \cl U$ and $\cl S^2\subset \cl V$. Then $S^{-1}S^{-1}\subset \cl V^{-1}$.
Lemma \ref{lgtp3} implies $\cl{S^{-1}}\subset \cl V^{-1}$. Since $H$ is dense in $G$ and $H$ is a topological group, $e\in\Int \cl{S^{-1}}$. Let $W=V\cap \Int \cl{S^{-1}}$. Then $W\in \nne$, $W^{-1}\subset \cl V$ and $W\subset V$. We get $WW^{-1}\subset \cl V^2\subset \cl U$.
\end{proof}

\begin{lemma}\label{lgtp4}
Let $G$ be a \rr-semitopological \goo{sqpe}-topological group, $U\in\tps$ and $M\subset G=\cl M$.
For any positive $\nom$ there are $x\in M$ and $V\in\nne$, so $V^nx\subset U$.
\end{lemma}
\begin{proof}
Let us prove it by induction on $n$.

$n=1$. Let $x\in U\cap M$ and $V=Ux^{-1}$.

$n>1$. By the induction hypothesis, there are $x_*\in U\cap M$ and $V_*\in \nne$ for which $V_*^{n-1}x_*\subset U$. Let $h\in V_*x_*$. Since $G$ is a \goo{sqpe}-topological group, there are $g\in V_*x_*h^{-1}$ and $S\in\nne$such that $S^2g\subset V_ *x_*h^{-1}$. Let $x\in Sgh\cap M$ and
\[
V = Sghx^{-1}\cap S\cap V_*.
\]
Then $V^2x\subset SSgh\subset V_*x_*$ and $V^nx\subset V_*^{n-1}x_*\subset U$.
\end{proof}

\begin{theorem}\label{tgtp5}
Let $G$ be a \rr-semitopological \goo{dfl^*,fpe}-topological $\gD$-Baire group.
\begin{itemize}
\item[{\rm (1)}]
If $G$ is a $\pi$-semiregular space, then $G$ is a $\gdd v$-Baire group;
\item[{\rm (2)}]
If $G$ is a semiregular space, then $G$ is a topological group.
\end{itemize}
\end{theorem}
\begin{proof}
The group $H=\tcs_f(G)$ is dense in $G$. Proposition \ref{pvc1} ($\msgoo{dfl^*,fpe,r}\to\msgoo{sqpe}$) implies that $G$ is a \goo{sqpe}-topological group, that is, multiplication in $G$ quasi-continuous in the first coordinate in $(e,e)$.

\begin{lemma}\label{lgtp5}
$\Int S^{-1}\neq\es$ for any regular open $S\in\tps$.
\end{lemma}
\begin{proof}
Lemma \ref{lgtp4} implies that there exists $g\in S\cap H$ and $U\in\nne$, so $U^2g\subset S$.
Proposition \ref{pgdp1} ($\gD$) implies that $WW^{-1}\subset \clg U$ for some $W\in\nne$. Since $\clg U\subset \cl{U^2}$, then
\[
WW^{-1}g\subset \clg U g \subset \cl{U^2} g \subset \cl S.
\]
Since $S$ is a regular open subset of $G$ and $WW^{-1}g$ is open, then $WW^{-1}g\subset S$. Then $g^{-1}WW^{-1}\subset S^{-1}$. Since $g\in H$, then $\Int g^{-1}WW^{-1}\neq \es$. Hence $\Int S^{-1}\neq\es$.
\end{proof}
Condition (2) follows from condition (1) and Theorem \ref{tgtp1}.
Let us prove (1).
Since $G$ is a $\pi$-semiregular space, it follows from the Lemma \ref{lgtp5} that $G$ is a \goo{fie}-topological group.
Proposition \ref{pvc1} ($\msgoo{dfl^*,fie,r}\to\msgoo{fl,fi}$) implies that $G$ is a \goo{fl,fi}-topological group.
So the group $G$ is a \goo{fl,sqpe,fi,r}-topological group. Let $U\in\nne$. Let us show that $\Int\clg U\subset\cl U$. Let us assume the opposite. Then $S=\Int\clg U\setminus\cl U\neq\es$.
Lemma \ref{lgtp4} implies that $q\in S$ and $Q\in\nne$ existsuch that $Q^2q\subset S$.
Since $G$ is a \goo{fi}-topological group, $\Int (Qq)^{-1}\neq\es$ and $Vg^{-1}\subset \Int (Qq)^{- 1}$ for some $g\in Qq$ and $V\in \nne$. We get $gV^{-1}\subset Qq$ and $QgV^{-1}\subset Q^2q\subset S$. Since $QgV^{-1}\cap U=\es$, then $Qg\cap UV=\es$ and $g\not\in \cl{UV}\supset \clg U$. A contradiction with the fact that $g\in S\subset \clg U$.

We have shown that $\Int\clg U\subset\cl U$. Proposition \ref{pgdp1} ($\gD$) implies that $WW^{-1}\subset \clg U$ for some $W\in\nne$. Since $WW^{-1}\in\tps$, then $WW^{-1}\subset \Int\clg U\subset\cl U$.
\end{proof}

\begin{proposition} \label{pgtp1}
Let $G$ be a \rr-semitopological \goo{sqpe}-topological $\gdd h$-Baire group. Then $G$ is a quasi-regular space.
\end{proposition}
\begin{proof} Let $U\in\tps$. Lemma \ref{lgtp4} implies that $V^3g\subset U$ for some $g\in U$ and $V\in\tps$. Since $\clh V\subset V^2$, it follows from Proposition \ref{pgdp1} ($\gdd h$) that $WW^{-1}\subset \clh V\subset V^2$ for some $W\in\nne$. Then $SVg\subset U$, where $S=WW^{-1}$. Since $S=S^{-1}$, it follows from the Lemma \ref{lgtp3} that $\cl V\subset SV$. Hence $\cl{Vg}\subset U$.
\end{proof}

\begin{theorem}\label{tgtp6}
Let $G$ be a \rr-semitopological \goo{dfl^*,fpe}-topological $\gdd h$-Baire group.
Then
\begin{itemize}
\item[{\rm (1)}]
$G$ is a $\gdd v$-Baire group;
\item[{\rm (2)}]
if $G$ is a semiregular space, then $G$ is a topological group.
\end{itemize}
\end{theorem}
\begin{proof}
Condition (1) follows from Theorem \ref{tgtp5} (1) and Proposition \ref{pgtp1}.
Condition (2) follows from condition (1) and Theorem \ref{tgtp1}. Also, condition (2) follows from Theorem \ref{tgtp5} (2).
\end{proof}

\begin{theorem}\label{tgtp7}
Let $G$ be a semitopological \goo{dl^*,fie}-topological $\gdd h$-Baire group.
Then
\begin{itemize}
\item[{\rm (1)}]
$G$ is a $\gdd v$-Baire group;
\item[{\rm (2)}]
if $G$ is a semiregular space, then $G$ is a topological group.
\end{itemize}
\end{theorem}
\begin{proof}
Condition (2) follows from condition (1) and Theorem \ref{tgtp1}. Let us prove (1).
Proposition \ref{pvc1} ($\msgoo{dfl^*,fie,r}\to\msgoo{fl,fi}$) implies that $G$ is a \goo{fi}-topological group.
The group $H=\tcs(G)$ is dense in $G$.
Let $U\in \nne$.
\begin{lemma}\label{lgtp6}
$\Int\clh U\subset \cl U$.
\end{lemma}
\begin{proof}
Let us assume the opposite. Then $S=\Int\clh U\setminus \cl U\neq\es$. Since $G$ is a \goo{fi}-topological group, $\Int S^{-1}\neq\es$. Then $Qg^{-1}\subset \Int S^{-1}$ for some $g\in H\cap S$ and $Q\in\nne$. Then $gQ^{-1}\subset S$ and $\cl{Q^{-1}}\subset g^{-1}\cl S$. Proposition \ref{pgdp1} ($\gdd h$) implies that $V=WW^{-1}\subset \cl{Q^{-1}}$ for some $W\in\nne$. Then $gV\subset \cl S$ and $gV\cap U=\es$. Hence $g\notin UV \supset \clh U$. A contradiction with $g\in S \subset \clh U$.
\end{proof}
Proposition \ref{pgdp1} ($\gdd h$) implies that $WW^{-1}\subset \clh U$ for some $W\in\nne$. Then $WW^{-1}\subset \Int\clh U\subset \cl U$.
\end{proof}

\section{Groups with quasicontinuous multiplication}\label{sec-CHART}


Let us define the topological games $G(y_*,Y)$ and $\tG(y_*,Y)$ for the space $Y$ and $y\in Y$ \cite{Gruenhage1976,DolezalMoors2017}.
Players $\al$ and $\be$ are playing. On the $n$th move, player $\al$ chooses
\begin{itemize}
\item open neighborhood $W_n\subset Y$ of point $y_*$ in game $G(y_*,Y)$;
\item is an open non-empty set $W_n\subset Y$ in the game $\tG(y_*,Y)$.
\end{itemize}
Player $\be$ chooses $y_n\in W_n$. Player $\al$ wins if $y_*\in\cl{\set{y_n:\nom}}$.

A point $x\in X$ is called a \term{$W$-point} (\term{$\tW$-point}) if if player $\al$ has a winning strategy in the game $G(y_*,Y )$ ($\tG(y_*,Y)$).
A space $X$ is called a $W$-space ($\tW$-space) if every point in $Y$ is a $W$-point ($\tW$-point) \cite{Gruenhage1976,DolezalMoors2017}.

Any $W$-space is a $\tW$-space. Property $W$ is inherited by subspaces, property $\tW$ is inherited by dense subspaces.

Recall the definition of the topological Banach-Mazur game $BM(X)$ for the space $X$. Denote $U_{-1}=X$.
Players $\al$ and $\be$ are playing. On the $n$th move, player $\be$ chooses an open non-empty set $V_n\subset U_{n-1}$, player $\al$ chooses an open non-empty set $U_n\subset V_{n}$. Player $\al$ wins if $\bigcap_\nom V_n=\bigcap_\nom U_n\neq\es$.

The Banach-Mazur theorem \cite{Oxtoby1957} says that if the player $\be$ has a winning strategy, then the space $X$ is not Baire.

A mapping $\Fp: X\times Y\to Z$ is called \term{KC-continuous} if the mappings $\Fp(\cdot,y_*): X\to Z,\ x\mapsto \Fp(x,y_ *)$ is quasicontinuous and the mappings $\Fp(x_*,\cdot): Y\to Z,\ y\mapsto \Fp(x_*,y)$ are continuous for all $(x_*,y_*)\in X\ times Y$.

\begin{proposition} \label{p:CHART:1}
Let $X$, $Y$, $Z$ be a regular space, the map $\Fp: X\times Y\to Z$ be KC-continuous, the point $y_*\in Y$ be a $\tW$-point, and the space $X$ Baire. Then the mapping $\Fp$ is quasi-continuous at every point $(x_*,y_*)$ for all $x_*\in X$.
\end{proposition}
\begin{proof} Assume the contrary. Then the mapping $\Fp$ is not quasi-continuous at the point $(x_*,y_*)$ for some point $x_*\in X$. There is a neighborhood $S\subset X\times Y$ of $(x_*,y_*)$ and a neighborhood $O\subset Z$ of $\Fp(x_*,y_*)$ so that $\Fp(S') \not \subset O$ for any non-empty $S'\subset S$.

Let $A$ and $A'$ be neighborhoods of the point $\Fp(x_*,y_*)$ such that $A\subset \cl A\subset A'\subset \cl{A'}\subset O$. Let $B= Z\setminus \cl{A'}$. Then $\cl A\cap \cl B=\es$ and $\Fp(S')\cap B\neq\es$ for any non-empty open $S'\subset S$. That is, $S\subset \cl{\Fp^{-1}(B)}$.
Take open non-empty open $U\subset X$ and $V\subset Y$ such that $y_*\in V$, $U\times V\subset S$ and $\Fp(U\times\sset{y_ *})\subset A$.

Let us construct a strategy $s^*_\be$ for player $\be$ in the Banach-Mazur game $BM(X)$. At the same time, we will play an auxiliary game in the game $\tG(y_*,Y)$, in which the player $\al$ follows the winning strategy $p_\al$.
We can assume that in the strategy $p_\al$ the player chooses an open subset $V$.
Let $s_\al$ be some strategy of player $\al$ in the game $BM(X)$. We will also construct an auxiliary strategy $p^*_\be$ of player $\be$ in the game $\tG(y_*,Y)$. The strategies $s_\al$ and $p_\al$ are given, the strategies $s^*_\be$ and $p^*_\be$ will be built.

During construction, the following conditions will be met:
\begin{itemize}
\item[($C_n$)] $\Fp(U_n\times \sset{y_n})\cap B\neq\es$;
\item[($D_n$)] $\Fp(V_{n+1}\times \sset{y_{n}})\subset B$.
\end{itemize}
{\bf $0$th step}:
\begin{itemize}
\item[($s^*_\be$)] $V_0=V$;
\item[($s_\al$)] $U_0\subset V_0$;
\item[($p_\al$)] $W_0\subset V$;
\item[($p^*_\be$)] $y_0\in W_0$ is chosen in such a way that the condition ($C_0$) is satisfied.
\end{itemize}
{\bf $n$th step}:
\begin{itemize}
\item[($s^*_\be$)] since the condition ($C_{n-1}$) is satisfied and the mapping $\Fp$ is KC-continuous, there exists a non-empty open $V_n\subset U_{n- 1}$ so that the condition ($D_{n-1}$) is satisfied;
\item[($s_\al$)] $U_n\subset V_n$;
\item[($p_\al$)] $W_n\subset V$;
\item[($p^*_\be$)] $y_n\in W_n$ is chosen such that ($C_n$) is satisfied.
\end{itemize}
Let's check that player $\be$ won the game $BM(X)$. Assume the opposite, that is, there exists $x\in \bigcap_\nom V_n$.
From ($D_n$) it follows that $\Fp(x,y_n)\in B$ for $\nom$.
Since the strategy $p_\al$ is winning, then $y_*\in\cl{\set{y_n:\nom}}$.
Since the function $\Fp(x,\cdot)$ is continuous, then $\Fp(x,y_*)\in \cl B$.
Since $\Fp(U\times\sset{y_*})\subset A$, then $\Fp(x,y_*)\in A$.
A contradiction with the fact that $\cl A\cap \cl B=\es$.

We have shown that the strategy $s^*_\be$ is winning for player $\be$ in the game $BM(X)$. Hence the space $X$ is not Baire. Contradiction.
\end{proof}

\begin{proposition} \label{p:CHART:2}
Let $G$ be a \goo{dql,r}-topological Baire regular group and $\tc_q(G)$ contain a dense $\tW$-space. Then $G$ is a \goo{qpe}-topological group.
\end{proposition}
\begin{proof}
Let $X$ be a $\tW$-space dense in $\tc_q(G)$. Since $G$ is homogeneous, we can assume that $e\in X$. Let's put
\[
\Fp: X\times G\to G, (x,g) \mapsto gx.
\]
The mapping $\Fp$ is KC-continuous. Proposition \ref{p:CHART:1} implies that $\Fp$ is quasi-continuous. Since $X$ is dense in $G$, the multiplication $\gm$ is quasicontinuous in $(e,e)$.
\end{proof}

Proposition \ref{p:CHART:2} implies the following statement.

\begin{proposition} \label{p:CHART:3}
Let $G$ be a \goo{ql,r}-topological Baire regular group and $G$ contain a dense $\tW$-space. Then $G$ is a \goo{qpe}-topological group.
\end{proposition}

\begin{proposition} \label{p:CHART:4}
A regular space with a countable $\pi$-character is a $\tW$-space.
\end{proposition}
\begin{proof}
Let $X$ be a space with countable $\pi$-character, $x_*\in X$, $\sqnn W$ be a countable $\pi$-base in $x_*$. The winning strategy for player $\al$ in the game $\tG(x_*,X)$ is that on the $n$th move player $\al$ chooses $W_n$ from the $\pi$-base.
\end{proof}

\begin{proposition} \label{p:CHART:5}
Let $G$ be a \goo{dql,r}-topological Baire regular group with countable $\pi$-character. Then $G$ is a \goo{qpe}-topological group.
\end{proposition}
\begin{proof}
The countable $\pi$-character is inherited by dense subspaces. Hence $\tc_q(G)$ has a countable $\pi$-character. Next, apply Proposition \ref{p:CHART:2}.
\end{proof}

\section{Main results}\label{sec-gmain}

In this section, we formulate corollaries from the Section \ref{sec-rtgc}.


Let $G$ be an right semitopological group.
A set $P\subset G\times G$ is called a semi-neighborhood of the diagonal if $P(x)=\set{y\in G: (x,y)\in P}$ is a neighborhood (not necessarily open) of a point $x\ in G$.
For a neighborhood $U\in \nne$ of unity, $\nn(U)=\set{(x,y)\in G^2: yx^{-1} \in U }$ is a semi-neighborhood of the diagonal.

For $\wt\D\in\sset{\D,\D_g,\D_s}$, a space $G$ (right semitopological group $G$) is $\wt\D$-nonmeager ($g\wt \D$-Baire) if for any semi-neighbourhood of the diagonal $P$
(semi-neighborhoods of the diagonal of the form $P=\nn(U)$, where $U$ is a neighborhood of unity) there exists an open non-empty $W\subset G$, so that the condition ($\wt\D$) is satisfied:
\begin{itemize}
\item[($\D$)] $W\times W\subset \cl{P\cap(W\times W)}$;
\item[($\D_h$)] $W\subset \cl{\set{x: (x,y)\in P}}$ for all $y\in W$;
\item[($\D_s$)] $W\subset \cl{\set{x: \sset x \times W \subset P}}$.
\end{itemize}

Let $\cP_k$ ($\cP_c$) be the smallest class of spaces that
\begin{itemize}
\item contains $p$-spaces and strongly $\Sigma$-spaces;
\item is closed under arbitrary (countable) products;
\item is closed under taking open subspaces.
\end{itemize}
Let $\cD_d$ be the smallest class of spaces that
\begin{itemize}
\item contains $\Sigma$-spaces, $w\D$-spaces, and feebly compact spaces;
\item is closed under products by spaces from the class $\cP_k$;
\item is closed under taking open subspaces.
\end{itemize}
Let $\cD_h$ be the smallest class of spaces that
\begin{itemize}
\item contains $\Sigma$-spaces and $w\D$-spaces;
\item is closed under products by spaces from the class $\cP_c$;
\item is closed under taking open subspaces.
\end{itemize}
Let $\cD_s$ be the class of semiregular spaces that contain a metrizable nonmeager subset. Baire semiregular spaces from the following classes belong to this class:
\begin{itemize}
\item $\sigma$-spaces, spaces with a countable network;
\item developable spaces.
\end{itemize}
In \cite{rezn2022gbd} $\godt_{o,l}$-nonmeager, $\godt_{p,l}$-nonmeager and $\gbmt_{f}$-nonmeager spaces are introduced and studied.

\begin{theorem} \label{tgmain-delta}
Let $G$ be a \rr-topological group.
\begin{itemize}
\item[{\rm ($\gD$)}]  
 A right semitopological group is $\gD$-Baire if it is a $\D$-nonmeager space. The class of $\D$-nonmeager spaces contains $\godt_{o,l}$-nonmeager spaces, which contain regular Baire spaces from the class $\cD_d$.
\item[{\rm ($\gD_h$)}] 
 A right semitopological group is $\gdd h$-Baire if it is a $\dd h$-nonmeager space. The class of $\dd h$-nonmeager spaces contains $\godt_{p,l}$-nonmeager spaces, which contain regular Baire spaces from the class $\cD_h$.
\item[{\rm ($\gD_s$)}] 
 A right semitopological group is $\gdd s$-Baire if it is $\dd s$-nonmeager or \cdp/-nonmeager space. The class of $\dd s$-nonmeager spaces contains $\gbmt_{f}$-nonmeager spaces. The classes of $\dd s$-nonmeager and \cdp/-nonmeager spaces contain the class $\cD_s$.
\end{itemize} 
\end{theorem}
\begin{proof} 
($\gD$)
It follows from Proposition \ref{pgdp1} that an right semitopological group is $\gD$-Baire if $G$ is a $\D$-nonmeager space. Proposition \dref{pgd1} (1) \cite{rezn2022gbd} implies that $\godt_{o,l}$-nonmeager spaces are $\D$-nonmeager. It follows from the Theorem \dref{tgd1} \cite{rezn2022gbd} that the class of $\godt_{o,l}$-nonmeager spaces contains regular Baire spaces from the class $\cD_d$.

($\gD_h$)
It follows from Proposition \ref{pgdp1} that an right semitopological group is $\gdd h$-Baire if $G$ is a $\dd h$-nonmeager space. Proposition \dref{pgd1} (2) \cite{rezn2022gbd} implies that $\godt_{p,l}$-nonmeager space is $\dd h$-nonmeager. It follows from the Theorem \dref{tgd1} \cite{rezn2022gbd} that the class of $\godt_{p,l}$-nonmeager spaces contains regular Baire spaces from the class $\cD_h$.

($\gD_s$)
It follows from Proposition \ref{pgdp1} that an right semitopological group is $\gdd s$-Baire if $G$ is a $\dd s$-nonmeager space. Proposition \dref{pgd1} (3) \cite{rezn2022gbd} implies that $\gbmt_{f}$-nonmeager space is $\dd s$-nonmeager. It follows from the Theorem \dref{tgd1} \cite{rezn2022gbd} that the class of $\gbmt_{f}$-nonmeager spaces contains regular Baire spaces from the class $\cD_s$. Proposition \ref{pdsp1} implies that \cdp/-nonmeager spaces are $\dd s$-nonmeager spaces. Proposition \ref{pcdp5} implies that the spaces in the class $\cD_s$ are \cdp/-nonmeager.
\end{proof}

\begin{theorem}[Theorem \ref{trtg2}]\label{tgmain0}
Let $G$ be a \rr-topological group. If the set $\tcs_f(G)$
dense in $G$, then $G$ is a topological group.
\end{theorem}

\begin{lemma} \label{lgmain1}
Let $G$ be a compact Hausdoff right semitopological group. Then $\tc(G)=\tcs(G)$.
\end{lemma}
\begin{proof}
Clearly, $\tc(G)\supset \tcs(G)$. Let $g\in \tc(G)$. The mapping $\la_g: G\to G$ is a bijective continuous mapping of compact Hausdoff spaces. Hence $\la_g$ is a homeormorphism and the mapping $(\la_g)^{-1}=\la_{g^{-1}}$ is continuous, i.e. $g^{-1}\in\tcs(G)$.
\end{proof}

\begin{theorem} \label{tgmain1}
Let $G$ be a $\gD$-Baire right semitopological group.
\begin{itemize}
\item[{\rm (1)}]
If $G$ is an \rr-paratopological group, then $G$ is an \rr-topological group.
\item[{\rm (2)}]
If $G$ is a semiregular space, the set $\tcs_f(G)$
dense in $G$ and
multiplication $\gm: G\times G\to G,\ (g,h)\mapsto gh$ is feebly continuous in $(e,e)$, then
 $G$ is a topological group.
\item[{\rm (3)}]
If $G$ is a regular space, the set $\tcs_q(G)$
is dense in $G$ and $\tcs_q(G)$ contains a dense $\tW$-space, then $G$ is a topological group.
\end{itemize}
\end{theorem}
\begin{proof}
Condition (1) follows from Theorem \ref{tgtp1} (5). Condition (2) follows from Theorem \ref{tgtp5}. Condition (3) follows from (2) and Proposition \ref{p:CHART:2}.
\end{proof}

\begin{cor} \label{cgmain1}
Let $G$ be a regular right semitopological group.
\begin{itemize}
\item[{\rm (1)}]
If $G$ is a product of \v Cech complete first-countable spaces (for example, $G$ is homeomorphic to $\R^\tau$) and $\la_g$ is quasi-continuous for every $g\in G$, then $G$ is a topological group.
\item[{\rm (2)}] \cite[Proposition 3.2]{Moors2016}
If $G$ is a CHART group with feebly continuous multiplication, then $G$ is a topological group.
\item[{\rm (3)}]
If $G$ is a CHART group with countable $\pi$-character (for example, $G$ is a compact space with countable tightness), then $G$ is a metrizable topological group.
\end{itemize}
\end{cor}
\begin{proof}
(1) The product of \v Cech complete spaces is a Baire space, and \v Cech complete spaces are $p$-spaces.
Theorem \ref{tgmain-delta} implies that $G$ is a $\gD$-Baire group. The product $G$ contains the $\Sigma$-product $S$ of first-countable spaces.
\cite[Theorem 4.6]{Gruenhage1976} implies that $S$ is a $W$-space and hence a $\tW$-space. Theorem \ref{tgmain1} (3) implies that $G$ is a topological group.

(2) Lemma \ref{lgmain1} implies that $\tc(G)=\tcs(G)\subset \tcs_f(G)$ and $\tcs_f(G)$ is dense in $G$. Theorem \ref{tgmain1} (2) implies that $G$ is a topological group.

(3) Compact spaces with countable tightness have countable $\pi$-character \cite{Shapirovskii1975}.
It follows from Proposition \ref{p:CHART:5} that $G$ is a \goo{qpe}-topological group and, therefore, multiplication in $G$ is feebly continuous. It follows from (2) that $G$ is a topological group.
Topological groups with countable $\pi$-character are first-countable \cite[Proposition 5.2.6]{at2009}. The Birkhoff--Kakutani theorem \cite[Theorem 3.3.12]{at2009} implies that $G$ is metrizable.
\end{proof}

\begin{theorem} \label{tgmain2}
Let $G$ be a semiregular $\gdd h$-Baire right semitopological group.
\begin{itemize}
\item[{\rm (1)}]
If $G$ is an \rr-quasitopological group, then
 $G$ is an \rr-topological group.
\item[{\rm (2)}]
If the set $\tcs(G)$
is dense in $G$ and the inverse mapping $\gi: G\to G,\ g\mapsto g^{-1}$ is feebly continuous in $e$, then
$G$ is a topological group.
\end{itemize}
\end{theorem}
\begin{proof}
Condition (1) follows from the Theorem \ref{tgtp1} (2) and (4). Condition (2) follows from Theorem \ref{tgtp7}.
\end{proof}

\begin{theorem} \label{tgmain3}
Let $G$ be an $\gdd s$-Baire right semitopological group.
\begin{itemize}
\item[{\rm (1)}]
If $\la_g$ is feebly continuous for any $g\in G$, then
 $G$ is a topological group.
\item[{\rm (2)}]
If $G$ is a semiregular space and the set $\tcs_f(G)$
dense in $G$, then
 $G$ is a topological group.
\end{itemize}
\end{theorem}
\begin{proof}
 Condition (1) follows from Theorem \ref{tgtp3}. Condition (2) follows from Theorem \ref{tgtp4}.
\end{proof}

\begin{cor} \label{cgmain2} \rm{(MA)}
Let $G$ be a CHART group and $w(G)<2^\omega$. Then $G$ is a topological group.
\end{cor}
\begin{proof}
CHART groups have a right-invariant Haar measure \cite{MilnesPym1992,Moors2015}. Hence $G$ is a ccc space. Corollary \ref{ccdp1} implies that $G$ is a \cdp/-space. Hence $X$ is a \cdp/-nonmeager space. Theorem \ref{tgmain-delta} and \ref{tgmain3} (2) implies that $G$ is a topological group.
\end{proof}

\section{Examples and questions}\label{sec-qe}

\begin{example}[Example \ref{ertg1}]\label{eqe1}
Metrizable compact \rr-topological non-topological group.
\end{example}

\begin{example}\label{eqe2}
Let $G=\R$ with a topology whose open sets are $U\setminus P$, where $U$ is open in $\R$ and $P$ is nowhere dense in $\R$. The group $G$ is a quasitopological \goo{sqpe}-topological quasi-regular Hausdorff $\gdd v$-Baire $\gdd h$-Baire $\D$-Baire non-$\dd h$-nonmeager group which is not a topological group.
\end{example}

\begin{example}\label{eqe3}
Let $G=\R^2$ with a topology whose base point $(x_*,y_*)\in G$ form sets of the form
\[
(U\cap \set{(x,y)\in G: y>y_*})\cup \sset{(x_*,y_*)},
\]
where $U$ is an open neighborhood $(x_*,y_*)$ in the standard topology $\R^2$ \cite{cdp2010}. The group $G$ is a Hausdorff quasi-regular non-regular paratopological group.
\end{example}

\begin{example}[Example 2.13 \cite{cdp2010}]\label{eqe4}
There exists a Baire and metric left semitopological group $G$ which
is not a right semitopological group, and whose inversion is nearly continuous but not
feebly continuous. The group $G$ is $\dd s$-Baire.
\end{example}

\begin{example}\label{eqe5}
Let $G$ be a pseudocompact Boolean quasitopological group that is not a topological group (Example \dref{eqe8} \cite{rezn2022gbd}). The group $G$ is a Tychonoff $\D$-Baire non-$\gdd h$-Baire group.
\end{example}

\begin{problem} \label{pqe1} Let $G$ be a (semi)regular \rr-semitopological group. Which of the following conditions imply that $G$ is an \rr-topological group?
\begin{itemize}
\item[{\rm (1)}]
The group $G$ is $\gD$-Baire ($\D$-Baire, pseudocompact) and the multiplication of $\gm$ in $G$ is quasi (feebly) continuous (in $(e,e)$).
\item[{\rm (2)}]
The group $G$ is $\gdd h$-Baire ($\dd h$-Baire, countably compact, compact) and the multiplication of $\gm$ in $G$ is quasi (feebly) continuous (in $(e,e)$) .
\item[{\rm (3)}]
The group $G$ is $\gdd s$-Baire ($\dd s$-Baire, metrizable Baire) and the multiplication of $\gm$ in $G$ is quasi (feebly) continuous (in $(e,e)$).
\item[{\rm (4)}]
The group $G$ is $\gdd h$-Baire ($\dd h$-Baire, countably compact, compact) and the operation of taking the inverse of $\gi$ in $G$ is quasi (feebly) continuous (in $e$).
\item[{\rm (5)}]
The group $G$ is $\gdd s$-Baire ($\dd s$-Baire, metrizable Baire) and the operation of taking the inverse of $\gi$ in $G$ is quasi (feebly) continuous (in $e$).
\end{itemize}
\end{problem}

\begin{problem} \label{pqe1+1} Let $G$ be a (semi)regular $\gdd h$-Baire ($\dd h$-Baire, countably compact, compact) and  \rr-semitopological group.
Which of the following conditions imply that $G$ is a topological group?
\begin{itemize}
\item[{\rm (1)}]
$G=\tc_f(G)$.
\item[{\rm (q)}]
$G=\tc_q(G)$.
\end{itemize}
\end{problem}

A space $X$ is called \term{weakly pseudocompact} if there exists a compact Hausdorff extension $bX$ of the space $X$ in which the space $X$ is $G_\de$-dense, i.e. $X$ intersects any nonempty $G_\de$ subset of $bX$ \cite{arh-rezn2005}. It is clear that the product of weakly pseudocompact spaces is weakly pseudocompact; in particular, the product of pseudocompact spaces is weakly pseudocompact.

\begin{problem} \label{pqe2} Let $G$ be a (semi)regular semitopological group. Which of the following conditions imply that $G$ is a topological group?
\begin{itemize}
\item[{\rm (1)}]
The group $G$ is $\gdd h$-Baire ($\dd h$-Baire).
\item[{\rm (2)}]
The group $G$ is paratopological and $G$ is weakly pseudocompact (product of pseudocompact spaces, product of two pseudocompact spaces) (Question \dref{pqe5} \cite{rezn2022gbd}).
\item[{\rm (3)}]
The group $G$ is feebly compact and belongs to one of the following classes of spaces: separable; countable tightness; $k$-space (Problem 3.5 \cite{tk2014}).
\end{itemize}
\end{problem}

\bibliographystyle{elsarticle-num}
\bibliography{rtg}

\begin{thebibliography}{10}
\expandafter\ifx\csname url\endcsname\relax
  \def\url#1{\texttt{#1}}\fi
\expandafter\ifx\csname urlprefix\endcsname\relax\def\urlprefix{URL }\fi
\expandafter\ifx\csname href\endcsname\relax
  \def\href#1#2{#2} \def\path#1{#1}\fi

\bibitem{Montgomery1936}
D.~{Montgomery}, {Continuity in topological groups}, {Bull. Am. Math. Soc.} 42
  (1936) 879--882.
\newblock \href {https://doi.org/10.1090/S0002-9904-1936-06456-6}
  {\path{doi:10.1090/S0002-9904-1936-06456-6}}.

\bibitem{Ellis1957}
R.~{Ellis}, {A note on the continuity of the inverse}, {Proc. Am. Math. Soc.} 8
  (1957) 372--373.
\newblock \href {https://doi.org/10.2307/2033747} {\path{doi:10.2307/2033747}}.

\bibitem{tk2014}
M.~Tkachenko, Paratopological and semitopological groups vs topological groups,
  Recent Progress in General Topology 3 (2014) 825--872.

\bibitem{Reznichenko1994}
E.~Reznichenko,
  \href{https://www.sciencedirect.com/science/article/pii/0166864194900213}{Extension
  of functions defined on products of pseudocompact spaces and continuity of
  the inverse in pseudocompact groups}, Topology and its Applications 59~(3)
  (1994) 233--244.
\newblock \href {https://doi.org/https://doi.org/10.1016/0166-8641(94)90021-3}
  {\path{doi:https://doi.org/10.1016/0166-8641(94)90021-3}}.
\newline\urlprefix\url{https://www.sciencedirect.com/science/article/pii/0166864194900213}

\bibitem{SoleckiSrivastava1997}
S.~Solecki, S.~Srivastava,
  \href{https://www.sciencedirect.com/science/article/pii/S0166864196001198}{Automatic
  continuity of group operations}, Topology and its Applications 77~(1) (1997)
  65--75.
\newblock \href {https://doi.org/https://doi.org/10.1016/S0166-8641(96)00119-8}
  {\path{doi:https://doi.org/10.1016/S0166-8641(96)00119-8}}.
\newline\urlprefix\url{https://www.sciencedirect.com/science/article/pii/S0166864196001198}

\bibitem{FerriHernandezWu2006}
S.~Ferri, S.~Hernández, T.~Wu,
  \href{https://www.sciencedirect.com/science/article/pii/S0166864105001318}{Continuity
  in topological groups}, Topology and its Applications 153~(9) (2006)
  1451--1457.
\newblock \href {https://doi.org/https://doi.org/10.1016/j.topol.2005.04.007}
  {\path{doi:https://doi.org/10.1016/j.topol.2005.04.007}}.
\newline\urlprefix\url{https://www.sciencedirect.com/science/article/pii/S0166864105001318}

\bibitem{cdp2010}
J.~Cao, R.~Drozdowski, Z.~Piotrowski, \href{http://eudml.org/doc/37996}{Weak
  continuity properties of topologized groups}, Czechoslovak Mathematical
  Journal 60~(1) (2010) 133--148.
\newline\urlprefix\url{http://eudml.org/doc/37996}

\bibitem{GlasnerMegrelishvili2013}
E.~Glasner, M.~Megrelishvili, Banach representations and affine
  compactifications of dynamical systems, in: Asymptotic geometric analysis,
  Springer, 2013, pp. 75--144.

\bibitem{Moors2016}
W.~B. Moors, Fragmentable mappings and chart groups, Fundamenta Mathematicae
  234 (2016) 191--200.

\bibitem{ruppert1975}
W.~Ruppert, Uber kompakte rechtstopologische gruppen mit gleichgradig stetigen
  linkstranslationen, Sitz. ber. d. Osterr. Akad. d. Wiss. Math.-naturw. Kl.
  184 (1975) 159--169.

\bibitem{namioka1972}
I.~Namioka, Right topological groups, distal flows, and a fixed-point theorem,
  Mathematical systems theory 6~(1) (1972) 193--209.

\bibitem{MoorsNamioka2013}
W.~B. Moors, I.~Namioka, Furstenberg’s structure theorem via chart groups,
  Ergodic Theory and Dynamical Systems 33~(3) (2013) 954--968.

\bibitem{at2009}
A.~Arhangel'skii, M.~Tkachenko,
  \href{https://doi.org/10.2991/978-94-91216-35-0}{Topological Groups and
  Related Structures}, Atlantis Press, 2008.
\newblock \href {https://doi.org/10.2991/978-94-91216-35-0}
  {\path{doi:10.2991/978-94-91216-35-0}}.
\newline\urlprefix\url{https://doi.org/10.2991/978-94-91216-35-0}

\bibitem{rezn2022gbtg}
E.~Reznichenko, Generalization of baire spaces using games, arXiv preprint
  arXiv:2203.02937 (2022) 35\href
  {https://doi.org/https://doi.org/10.48550/arXiv.2203.02937}
  {\path{doi:https://doi.org/10.48550/arXiv.2203.02937}}.

\bibitem{rezn2022gbd}
E.~{Reznichenko}, {Generalization of Bair spaces using diagonal}, arXiv
  e-prints (2022) 35\href {http://arxiv.org/abs/2203.09389}
  {\path{arXiv:2203.09389}}, \href
  {https://doi.org/https://doi.org/10.48550/arXiv.2203.09389}
  {\path{doi:https://doi.org/10.48550/arXiv.2203.09389}}.

\bibitem{Ravsky2001}
O.~Ravsky, Paratopological groups ii, Matematychni Studii 16~(1) (2001)
  93--101.

\bibitem{Ameen2021}
Z.~A. Ameen, \href{https://doi.org/10.2478/mjpaa-2021-0009}{Almost somewhat
  near continuity and near regularity}, Moroccan Journal of Pure and Applied
  Analysis 7~(1) (2021) 88--99.
\newblock \href {https://doi.org/doi:10.2478/mjpaa-2021-0009}
  {\path{doi:doi:10.2478/mjpaa-2021-0009}}.
\newline\urlprefix\url{https://doi.org/10.2478/mjpaa-2021-0009}

\bibitem{Piotrowski1980}
Z.~Piotrowski, Quasi-continuity and product spaces, Geometric topology, {Proc}.
  int. {Conf}., {Warszawa} 1978, 349-352 (1980). (1980).

\bibitem{Bouziad1996}
A.~Bouziad, Continuity of separately continuous group actions in $p$-spaces,
  Topology Appl. 71~(2) (1996) 119--124.
\newblock \href {https://doi.org/10.1016/0166-8641(95)00039-9}
  {\path{doi:10.1016/0166-8641(95)00039-9}}.

\bibitem{KenderovKortezovMoors2001}
P.~S. Kenderov, I.~S. Kortezov, W.~B. Moors, Topological games and topological
  groups, Topology Appl. 109~(2) (2001) 157--165.
\newblock \href {https://doi.org/10.1016/S0166-8641(99)00152-2}
  {\path{doi:10.1016/S0166-8641(99)00152-2}}.

\bibitem{arh-rezn2005}
A.~Arhangel'skii, E.~Reznichenko, Paratopological and semitopological groups
  versus topological groups, Topology and its Applications 151~(1-3) (2005)
  107--119.

\bibitem{BinOst2008}
N.~Bingham, A.~Ostaszewski, Normed versus topological groups: Dichotomy and
  duality, Dissertationes Mathematicae 472 (2008) 151.
\newblock \href {https://doi.org/10.4064/dm472-0-1}
  {\path{doi:10.4064/dm472-0-1}}.

\bibitem{milnes1993}
P.~Milnes, Representations of compact right topological groups, Canadian
  Mathematical Bulletin 36~(3) (1993) 314–323.
\newblock \href {https://doi.org/10.4153/CMB-1993-044-1}
  {\path{doi:10.4153/CMB-1993-044-1}}.

\bibitem{rezn2008}
E.~Reznichenko, \href{http://gtopology.math.msu.su/erezn-ci}{Continuity of the
  inverse}, in russian (2008).
\newblock \href {http://arxiv.org/abs/2106.01803} {\path{arXiv:2106.01803}}.
\newline\urlprefix\url{http://gtopology.math.msu.su/erezn-ci}

\bibitem{mtw2007}
J.~{van Mill}, V.~Tkachuk, R.~Wilson,
  \href{https://www.sciencedirect.com/science/article/pii/S0166864107000193}{Classes
  defined by stars and neighbourhood assignments}, Topology and its
  Applications 154~(10) (2007) 2127--2134, special Issue: The 6th Iberoamerican
  Conference on Topology and its Applications (VI-CITA).
\newblock \href {https://doi.org/https://doi.org/10.1016/j.topol.2006.03.029}
  {\path{doi:https://doi.org/10.1016/j.topol.2006.03.029}}.
\newline\urlprefix\url{https://www.sciencedirect.com/science/article/pii/S0166864107000193}

\bibitem{moors2017}
W.~B. Moors,
  \href{https://www.sciencedirect.com/science/article/pii/S0166864117304273}{Some
  baire semitopological groups that are topological groups}, Topology and its
  Applications 230 (2017) 381--392.
\newblock \href {https://doi.org/https://doi.org/10.1016/j.topol.2017.08.042}
  {\path{doi:https://doi.org/10.1016/j.topol.2017.08.042}}.
\newline\urlprefix\url{https://www.sciencedirect.com/science/article/pii/S0166864117304273}

\bibitem{Gruenhage1976}
G.~Gruenhage,
  \href{https://www.sciencedirect.com/science/article/pii/0016660X76900246}{Infinite
  games and generalizations of first-countable spaces}, General Topology and
  its Applications 6~(3) (1976) 339--352.
\newblock \href {https://doi.org/https://doi.org/10.1016/0016-660X(76)90024-6}
  {\path{doi:https://doi.org/10.1016/0016-660X(76)90024-6}}.
\newline\urlprefix\url{https://www.sciencedirect.com/science/article/pii/0016660X76900246}

\bibitem{DolezalMoors2017}
M.~Doležal, W.~B. Moors,
  \href{https://www.sciencedirect.com/science/article/pii/S0166864117304406}{On
  a certain generalization of w-spaces}, Topology and its Applications 231
  (2017) 1--9.
\newblock \href {https://doi.org/https://doi.org/10.1016/j.topol.2017.09.001}
  {\path{doi:https://doi.org/10.1016/j.topol.2017.09.001}}.
\newline\urlprefix\url{https://www.sciencedirect.com/science/article/pii/S0166864117304406}

\bibitem{Oxtoby1957}
J.~C. Oxtoby, The Banach-Mazur game and Banach category theorem, Princeton
  University Press, Princeton, NJ, 1957, pp. 159--163.

\bibitem{Shapirovskii1975}
B.~{\`E}. Shapirovskii, On $\pi$-character and $\pi$-weight in compact
  hausdorff spaces, Doklady Akademii Nauk 223~(4) (1975) 799--802.

\bibitem{MilnesPym1992}
P.~Milnes, J.~Pym,
  \href{https://www.scopus.com/inward/record.uri?eid=2-s2.0-84968480767&doi=10.1090%2fS0002-9939-1992-1065088-1&partnerID=40&md5=8eba4cf3993aee61498e0845ba31cf4d}{Haar
  measure for compact right topological groups}, Proceedings of the American
  Mathematical Society 114~(2) (1992) 387 – 393, cited by: 15; All Open
  Access, Bronze Open Access, Green Open Access.
\newblock \href {https://doi.org/10.1090/S0002-9939-1992-1065088-1}
  {\path{doi:10.1090/S0002-9939-1992-1065088-1}}.
\newline\urlprefix\url{https://www.scopus.com/inward/record.uri?eid=2-s2.0-84968480767&doi=10.1090%2fS0002-9939-1992-1065088-1&partnerID=40&md5=8eba4cf3993aee61498e0845ba31cf4d}

\bibitem{Moors2015}
W.~B. Moors, Invariant means on chart groups, Khayyam Journal of Mathematics
  1~(1) (2015) 36--44.

\end{thebibliography}
\end{document}